\documentclass[12pt]{article}

\usepackage{textcomp, amsmath, amssymb, amsthm, enumerate}
\usepackage{graphicx}
\usepackage{yhmath}

\newtheorem{theorem}{Theorem}[section]
\newtheorem{lemma}[theorem]{Lemma}

\def\rr{{\mathbb R}}

\def\sik{{\rr}^2}

\def\su{\subset}

\def\se{\setminus}

\def\al{\alpha}
\def\be{\beta}
\def\ga{\gamma}

\def\de{\delta}
\def\De{\Delta}
\def\ep{\varepsilon}
\def\si{\sigma}
\def\Si{\Sigma}

\def\la{\lambda}
\def\La{\Lambda}

\def\cd{\cdot}
\def\stb{,\ldots ,}
\def\ha{\frac{1}{2}}
\def\emp{\emptyset}
\def\del{\partial}

\def\nl{[0,1]}

\def\msk{\medskip}
\def\bsk{\bigskip}
\def\noi{\noindent}

\def\ol{\overline}

\def\sumj0m{\sum_{j=0}^m}
\def\sumik{\sum_{i=1}^k}

\def\sumi0n{\sum_{i=0}^n}

\begin{tiny}\end{tiny}

\def\proof{\noi {\bf Proof.} }

\def\x1n{x_1 \stb x_n}
\def\y1n{y_1 \stb y_n}

\def\diam{{\rm diam}\, }
\def\dist{{\rm dist}\, }

\def\iv{\wideparen}

\date{August 12, 2016}

\begin{document}

\title{THE KAKEYA PROBLEM FOR CIRCULAR ARCS}

\author{K. H\'era and M. Laczkovich}

\maketitle

\footnotetext[3]{The research of the second author was partially
supported by the Hungarian
National Research, Development and Innovation Office,
Grant No. NKFIH 104178}

\begin{abstract}
We  prove that if a circular arc has angle short enough, 
then it can be continuously moved to any prescribed position
within a set of arbitrarily small area.
\end{abstract}

\section{Introduction and main results}

It is well-known that a line segment can be continuously moved in a plane
set of arbitrarily small area such 
that it returns to its starting position,
with its direction reversed. 
This fact was first proved by Besicovitch as a solution
to the classical Kakeya problem \cite{B}. Besicovitch's construction used 
the so-called `P\'al joins' in order to shift
the line segment to an arbitrary parallel position using 
arbitrarily small area. Using these P\'al joins, one can easily
deduce from Besicovitch's theorem that every line segment
can be moved to arrive at any prescribed position
within a set of arbitrarily small area.

In the paper `Three Kakeya problems' 
F. Cunningham asked whether or not a circular arc has a similar
property \cite[p. 591]{C2}. Our aim is to show that
the answer is affirmative, at least for circular arcs of angle short enough.

\begin{theorem}\label{t1}
Let $\si$ be a circular arc in the plane 
of radius $1$ and of arc length less than $1.32$. Then for every $\ep >0$,
we can move $\si$ continuously to any given 
position such that the area of the points touched by the
moving arc is less than $\ep$.
\end{theorem}

In \cite{HL} we introduced the following terminology.
We say that the set $A\su \sik$ has property (K), if there 
exist two different positions of $A$ such that $A$
can be continuously moved from the first position to
the second within a set of arbitrarily small area.
It is obvious that all line segments and circles have property (K).

The set $A$ is said to have property ${\rm (K^s )}$ if it 
can be moved to arrive at any prescribed position
within a set of arbitrarily small area.
As we saw above, line segments have property ${\rm (K^s )}$. According to
Theorem \ref{t1}, those circular arcs which have angle short enough
also have property ${\rm (K^s )}$. On the other hand, 
(full) circles do not have property ${\rm (K^s )}$ since every continuous movement
placing the circle far enough must touch every point inside the circle.
We do not know whether or not all circular arcs 
shorter than the full circle have property ${\rm (K^s )}$, as our
construction does not seem to work in the case where the angle of the arc is longer 
than a certain bound. We remark that 
apart from line segments and circular arcs no continuum can have
property ${\rm (K^s )}$ (see \cite[Theorem 1.2]{HL}).

We can easily reduce Theorem \ref{t1}
to the case when the final position is obtained from $\si$ by rotating it
about one of its endpoints.

\begin{theorem}\label{t2}
Let $\si _0$ be a circular arc of radius $1$ and of arc length 
less than $1.32$, and let $\si_1$ be obtained 
from $\si _0$ by rotating it about one of its  endpoints.
Then, for every $\ep >0$,
we can move $\si _0$ continuously to $\si _1$ such that the area of 
the points touched by the moving arc is less than $\ep$.
\end{theorem}

Assuming Theorem \ref{t2}, we can prove Theorem \ref{t1} as follows.
Let $\si _0$ and $\si _1$ be congruent circular arcs of radius
$1$ and of arc length less than $1.32$. Let $O_0$ and $O_1$
be the centres of the circles of $\si _0$ and $\si _1$, respectively.
Let $O_0 =P_0 ,P_1 \stb P_n =O_1$ be points such that
$0< |P_{i-1}P_i |<2$ for every $i=1\stb n$, where $|AB|$ denotes
the distance between the points $A$ and $B$.
Let $K_i$ denote the circle with centre $P_i$ and radius $1$
$(i=0\stb n)$. Then, for every $i=1\stb n$, the circles $K_{i-1}$ and $K_i$
have two common points; let $M_i$ be one of them.

Let $\ep >0$ be given. First we rotate the arc $\si$ with
initial position $\si _0$ about the point $O_0 =P_0$ until one of its
endpoints becomes $M_1$. Then, using Theorem \ref{t2}, we apply
a suitable continuous motion 
touching an area less than $\ep /n$ which places $\si$ into the circle
$K_1$. Then we rotate $\si$ about the point $P_1$ until one of its
endpoints becomes $M_2$. Using Theorem \ref{t2}, we apply
a suitable continuous motion 
touching an area less than $\ep /n$ which places $\si$ into the circle
$K_2$. Iterating this process $n$ times, $\si$ will be moved to a subarc
of $K_n$. Then, using a suitable rotation about the point $P_n =O_1$ 
we reach the final position $\si _1$.
It is clear that the area of the set of points touched by this motion
is less than $\ep$, which proves Theorem \ref{t1}.

We shall prove Theorem \ref{t2} in the next 5 sections.
After stating some preliminary lemmas in the next section, we
describe the basic construction, an
adaptation of Cunningham's sprouting process (see \cite[pp. 118-120]{C1}),
in Section \ref{s2}.
We estimate the area touched by the moving arc in Section \ref{s3},
and then, using this estimate, we prove Theorem \ref{t2} 
in Section \ref{s4}. The proofs of the lemmas stated and used during the
proof are given in Section \ref{s5}.

\section{Preliminary lemmas}

\begin{lemma}\label{l1}
Let $K$ be a circle of radius $1$ and centre $O$, 
and let $P$ be a point such that
$\dist (P,K) =d<1/2$. 
Let $Q\in K$ be a point such that 
the orientation of the triangle $\De_{POQ}$ is positive
(counter-clockwise). 
\begin{enumerate}[{\rm (i)}]
\item If $P$ is in the interior of $K$ and 
$|PQ|\ge 2\sqrt d$, then there is an angle $0<\al <\pi /2$
such that $\sin \al <2 \sqrt{d}$, and 
rotating $K$ about the point $Q$ in the positive direction by
angle $\al$, the circle $K'$ obtained contains $P$.
\item If $P$ is in the exterior of $K$ and 
$3\sqrt d \le |PQ|\le 2-d$, then there is an angle $0<\al <\pi /2$
such that $\sin \al <2 \sqrt{d}$, and 
rotating $K$ about the point $Q$ in the negative direction by
angle $\al$, the circle $K'$ obtained contains $P$.
\end{enumerate}
Moreover, the distance
between $O$ and the centre of $K'$ is less than $4\sqrt d$.
\end{lemma}

Figures \ref{fig:in} and \ref{fig:out} on pp. \pageref{fig:in} and 
\pageref{fig:out} illustrate the two cases of Lemma \ref{l1}.

\begin{lemma}\label{l3}
Let $K$ be a circle of radius $1$ and centre $O$. Let $A,B\in K$
be such that the orientation of the triangle $\De_{OAB}$ is positive,
and $\eta =BOA\angle <1/5$. Let $K_A$ and $K_B$ be the circles
obtained by rotating $K$ about the points $A$ and $B$ 
in the positive direction
by the angles $\al$ and $\be$,
respectively, where $0<\al <3\be /4$ and $\be <\eta$.
Then one of the intersection points of $K_A$ and $K_B$  
is inside $K$, and its distance from $A$ is less than $20\eta$.
\end{lemma}

See Figure \ref{fig:inout} on page \pageref{fig:inout}.

\begin{lemma}\label{l9}
Let $K$ be a circle of radius $1$ and centre $O$. Let $A,B\in K$
be such that the orientation of the triangle $\De_{OAB}$ is positive,
and $\eta =BOA\angle <1/10$. Let $K_A$ and $K_B$ be the circles
obtained by rotating $K$ about the points $A$ and $B$ 
in the negative direction
by the angles $\al$ and $\be$,
respectively, where $0<\al <3\be /4$ and $\be <\eta$.
Then one of the intersection points of $K_A$ and $K_B$  
is outside $K$, and its distance from $A$ is less than $50\eta$.
\end{lemma}

See Figure \ref{fig8} on page \pageref{fig8}.

\begin{lemma}\label{l6}
Let $K$ be a circle of radius $r$, and let $\iv{AB}$ be a subarc of $K$
of length $<r\pi$. Rotating the arc $\iv{AB}$ about the point $A$ by an
angle $\al$ we obtain the arc $\iv{AC}$. Then the area of the domain
$H$ bounded by the arcs $\iv{AB}$, $\iv{AC}$ and $\iv{BC}$ is 
$|AB|^2 \cd \al /2$.
\end{lemma}

\bsk
In the next two lemmas we shall use the 
following notation. Let $h , \ep$ be positive numbers satisfying
$\ep <10^{-6}$ and $h\le \ep /10^{3}$. 
Let $K_0$ and $K_1$ be the circles of radius $1$ and 
centres $(-\sin (h /2),0)$ and 
$(\sin (h/2),0)$, respectively. The point with coordinates
$(0, \cos (h/2))$ will be denoted by $M$. Thus $M\in K_0 \cap K_1$.

Let $\ol K _i$ denote the disc bounded by $K_i$ $(i=0,1)$.
The closures of the sets $\ol K _0 \se \ol K _1$ and 
$\ol K _1 \se \ol K _0$ are denoted by $L_0$ and $L_1$. 

The circle of radius $\ep$ and centre $M$ intersects
the lune $L_0$ in the arc $\iv{P_0  P_1}$, where $P_0 \in K_0$ and 
$P_1 \in K_1$.

\begin{lemma}\label{l4}
If a line intersects the arc $\iv{P_0  P_1}$ at two points, 
or touches the arc $\iv{P_0  P_1}$ ,
then the angle between the line and the $y$ axis is less than $h+\ep$.
\end{lemma}

\begin{lemma}\label{l5}
Let $Q'$ and $Q''$ be points of the arc $\iv{P_0  P_1}$
such that the $y$ coordinate of $Q'$ is greater than 
the $y$ coordinate of $Q''$. Let $C\in L_1$ such that 
$|CM|\le 2-5\ep$, and 
let $K'$ and $K''$ be circles with the following properties. The circle
$K'$ has radius $1$ and centre $O'$, contains the points $Q'$ and $C$,
the circle $K''$ has radius $1$ and centre $O''$, and contains the 
points $Q''$ and $C$, and $|OO'|, |OO''|<\ep$. 

Let $\rho$ denote the rotation about the point $C$ mapping $K''$ onto 
$K'$ and let the angle of the rotation $\rho$ be $\al$. Then 
\begin{enumerate}[{\rm (i)}]
\item $\al <\ep$,
\item the angle between the $x$ axis and the line going through the points $Q'$
and $\rho (Q'')$ is less than $6\ep$.
\end{enumerate} 
\end{lemma}

\section{The sprouting process} \label{s2}

First we fix three parameters: the positive numbers
$h , \ep$ satisfying $\ep <10^{-6}$ and $h<\ep /10^{3}$, 
and the positive integer $n$.

Let $K_0$ and $K_1$ be the circles of radius $1$ and 
centres $O_0 =(-\sin (h/2),0)$ and 
$O_1 =(\sin (h/2),0)$, respectively. The intersection $K_0 \cap K_1$ consists
of the points $M,N$ lying on the $y$ axis, where the $y$ coordinate
of $M$ is positive, and the $y$ coordinate
of $N$ is negative. Then the triangle with vertices $O_0 ,M, O_1$
is isosceles, and $O_0 M O_1 \angle =h$. This means that
rotating $K_0$ about the point $M$ by angle $h$ 
in the positive direction we obtain $K_1$.

We put $R=1.227 < 2-5\ep$. Since $|MN|=2\cos (h/2)>2-h>2-\ep$, we have
$R< |{MN}|$.

Let $\ol K _i$ denote the disc bounded by $K_i$ $(i=0,1)$.
The closures of the sets $\ol K _0 \se \ol K _1$ and 
$\ol K _1 \se \ol K _0$ are denoted by $L_0$ and $L_1$. 

The circle of radius $\ep$ and centre $M$ intersects
the lune $L_0$ in the arc $\iv{P_0  P_1}$, 
where $P_0 \in K_0$ and $P_1 \in K_1$.
The length of $\iv{P_0  P_1}$ equals $h\ep$, since $P_0 M P_1 \angle =h$.

For every $x\in \nl$, $P_x$ will denote the point of the 
arc $\iv{P_0 P_1}$ such that the length of the arc $\iv{P_0 P_x}$
equals $x$ times the length of the arc $\iv{P_0 P_1}$. Thus
$P_{1/2}$ is the middle point of the arc $\iv{P_0 P_1}$. 

We put $r_i =iR/n$ $(i=1\stb n)$.
The circle of radius $r_i$ and centre $M$ intersects
the lune $L_1$ in the arc $\iv{A_i^0 A_i^1}$, where $A_i^0 \in K_0$ and 
$A_i^1 \in K_1$ $(i=1\stb n)$. 
We define $A_0^0 =A_0^1 =M$.

We denote by $H$ the horn shaped domain
contained in $L_1$, bounded by the circles $K_0 , K_1$ and the circle 
of radius $R$ and centre $M$. Then the diameter of $H$ is the distance
$|{A_n^0 M}|=|{A_n^1 M}|=R$.

The `sprouting process' is a construction of a system of circles
of radius $1$ inductively. We start with the circles
$K_0$ and $K_1$. We put $K_0^0 =K_0$ and $K_1^1 =K_1$.
We rotate $K_0^0$ about the point $A^0_1$ such that
the resulting circle, denoted by $K_1^{1/2}$, contains $P_{1/2}$.
Similarly, we rotate $K_1^1$ about the point $A_1^1$ such that 
the resulting circle, denoted by $K_0^{1/2}$, contains $P_{1/2}$.

In this way we have defined the circles $K_0^0 ,K_1^{1/2} ,K_0^{1/2}, K_1^1$
such that $P_0 \in K_0^0$, $P_{1/2}\in K_0^{1/2} \cap K_1^{1/2}$
and $P_1 \in  K_1^1$. 
Using these circles we can transform $K_0$ to $K_1$ by rotating
$K_0 =K_0^0$ about $A^0_1$ to obtain $K_1^{1/2}$, then rotating the latter
about the point $P_{1/2}$ to obtain $K_0^{1/2}$, and then rotating the latter
about $A_1^1$ to obtain $K_1^1 =K_1$. 

\begin{figure}
\centering
\includegraphics[width=4in]{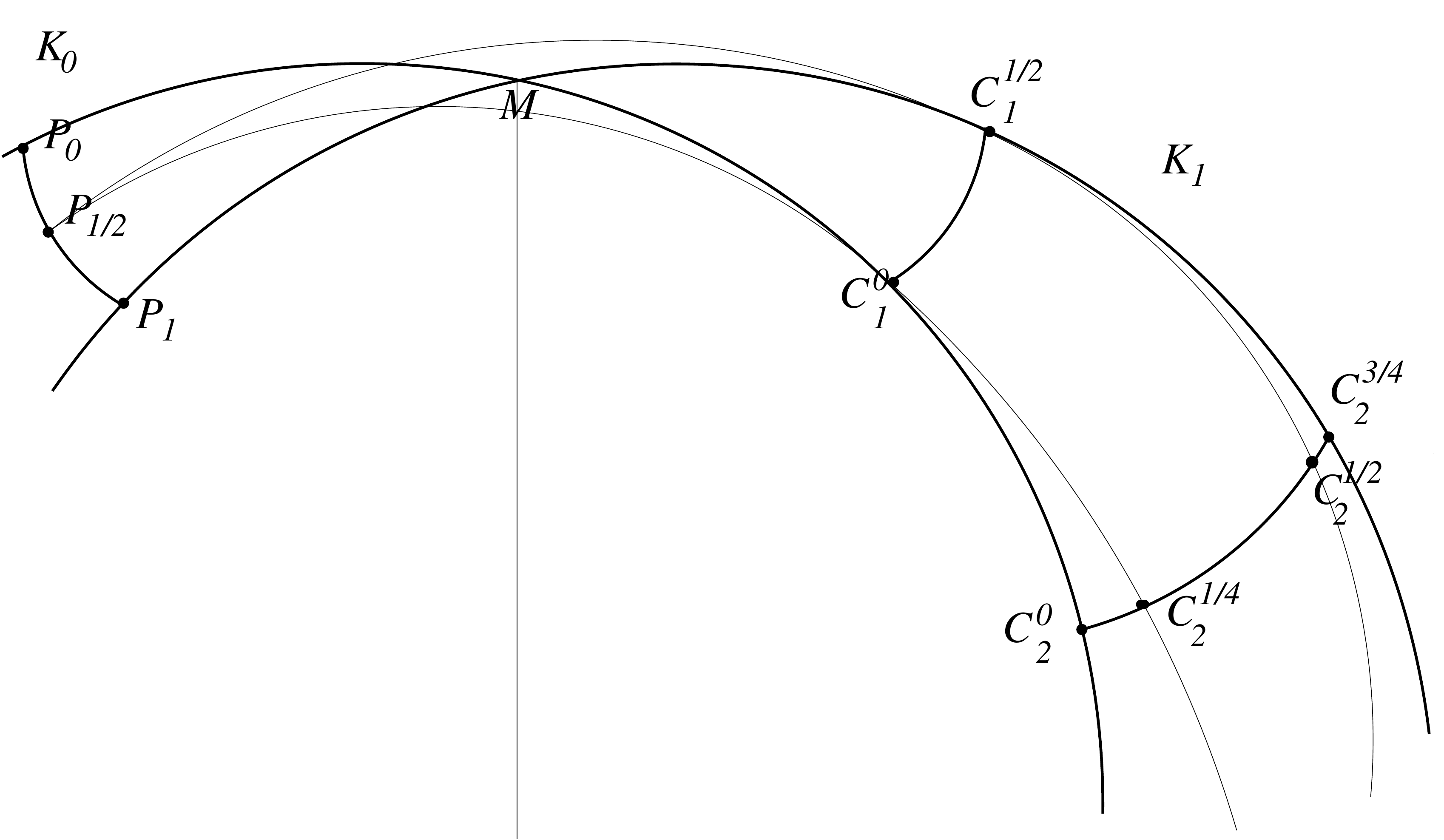}
\caption{Sprouting}
\label{fig12}
\end{figure}

We continue the process as follows. Let 
$K_0^0 ,K_1^{1/2} ,K_0^{1/2}, K_1^1$ intersect the arc $\iv{A^0_2 A^1_2}$
at the points $C_2^0 ,C_2^{1/4}, C_2^{1/2}, C_2^{3/4}$, respectively.
Thus $C_2^{0} =A^0_2$ and $C_2^{3/4} =A^1_2$ (see Figure \ref{fig12}).

Rotating $K_0^0$ about $C_2^0$ we obtain $K_1^{1/4}$ containing
$P_{1/4}$, rotating $K_1^{1/2}$ about $C_2^{1/4}$ we obtain $K_0^{1/4}$
containing $P_{1/4}$, rotating $K_0^{1/2}$ about $C_2^{1/2}$ we obtain $K_1^{3/4}$
containing $P_{3/4}$, 
and rotating $K_1^{1}$ about $C_2^{3/4}$ 
we obtain $K_0^{3/4}$ containing $P_{3/4}$. Then, using the circles
$K_0^0 \stb K_0^{3/4}$ and $K_1^{1/4} \stb K_1 ^1$ 
we can transform $K_0$ to $K_1$ by using rotations about the points
$C_2^0 \stb C_2^{3/4}$ and $P_{1/4}, P_{1/2}, P_{3/4}$ alternately.

Now we describe the general construction in detail.
We denote by $D_i$ the set $\{ k/2^i : 0\le k\le  2^i \}$ for every 
$i=0,1,\ldots$. We define the point $C_{i}^x$ and the circle $K^x_0$ 
for every $0\le i\le n$, $x\in D_i$, $x<1$, 
and the circle $K^x_1$ for every $0\le i\le n$, $x\in D_i$, $x>0$
satisfying the following conditions.
For every $0\le i\le n$ we have
\begin{equation}\label{e2}
C_i^x \in \iv{A_i^0 A_i^1} \qquad (x\in D_i , \ x<1 );
\end{equation}
\begin{equation}\label{e3}
K^{x}_0\ \text{contains the points}\  P_{x} \ \text{and}\ C_{i}^{x} 
\qquad (x\in D_i , \ x<1);
\end{equation}
\begin{equation}\label{e4}
K^x_1\ \text{contains the points}\  P_{x} \ \text{and}\ C_{i}^{x-2^{-i}} 
\qquad (x\in D_i , \ x>0);
\end{equation}
the centres of the circles $K_0^x$ and $K_1^{x+2^{-i}}$ 
are closer to the origin than
\begin{equation}\label{e5}
h+ 3\sqrt{h\ep} 
\cd \sum_{j=0}^{i-1} 2^{-j/2} \qquad (x\in D_i , \ x<1 ).
\end{equation} 
In addition,
\begin{equation}\label{e6}
\begin{split}
&K^x_0 \ \text{and} \ K_1^{x+2^{-i}}\ \text{intersect} \ \text{the closed arc} \\
&\qquad  \iv{A_j^0 A_j^1} \ (x\in D_i , \ x<1 , \ i\le j\le n).
\end{split}
\end{equation} 

We put $K^{0}_0 =K_0$ and $K^{1}_1 =K_1$, $C_0^0 =M$.
It is clear that \eqref{e2}-\eqref{e6} are satisfied for $i=0$.
(As for \eqref{e5}, we interpret the empty sum in \eqref{e5} as zero.
Note also that $|O_0 O|=|O_1 O|=\sin (h/2)<h$, where $O$ denotes the origin.)

Let $0\le i <n$ be given, and suppose that 
the points $C_i^x$ and the
circles $K^{x}_0 , K^{x}_1$ have been defined 
for the relevant values of $x$,
and satisfy \eqref{e2}-\eqref{e6}.
 
For every $x\in D_i$, $x<1$
we denote the intersection of $K^{x}_0$ and $\iv{A_{i+1}^0 A_{i+1}^1}$ by
$C_{i+1}^x$. Similarly, for every $x\in D_i$, $x>0$
we denote the intersection of $K^{x}_1$ and $\iv{A_{i+1}^0 A_{i+1}^1}$ by
$C_{i+1}^{x-2^{-i-1}}$. These points exist by \eqref{e6}.

Let $x\in D_i$, $x>0$ be fixed.
The point $Q=P_{x-2^{-i-1}}$ is the middle point
of the arc $\iv{P_{x-2^{-i}} P_x}$.
Therefore, $Q$ is outside the circle $K^x_1$. 
The distance $d$ between the point $Q$ and the circle 
$K^x_1$ is less than the length
of the arc $\iv{QP_x}$, which is $h\ep  /2^{i+1}$.
The distance $d'$ between $C_{i+1}^{x-2^{-i-1}}$ and $Q$
is greater than $\ep >3\sqrt{h\ep} >3 \sqrt d$ by $h<\ep /10$.
In addition, we have $d'\le R+\ep <2-\ep < 2-d$.

Therefore, we may apply
(ii) of Lemma \ref{l1}, and find that the circle 
$K^x_1$ can be rotated about the point $C_{i+1}^{x-2^{-i-1}}$ by an angle $\al$
in the negative direction such that $\sin \al <2\sqrt d$
and the rotated copy
of $K^x_1$, denoted by $K^{x-2^{-i-1}}_0$, contains the point $Q$.
By the lemma, the distance between the centres of $K^x_1$ and 
$K^{x-2^{-i-1}}_0$ is less than 
\begin{equation*}
4\sqrt d <\frac{4}{\sqrt 2} \cd \sqrt{h\ep} \cd 2^{-i/2} < 3\sqrt{h\ep} 
\cd 2^{-i/2}.
\end{equation*} 
Consequently, 
the distance between the origin and the centre of
$K^{x-2^{-i-1}}_0$ is less than $h +3\sqrt{h\ep}  \cd \sum_{j=0}^i 2^{-j/2}$.

In this way we have defined the circles $K_0^y$ for every 
$y\in D_{i+1} \se D_i$.

In the same way, we obtain the circle $K^{x+2^{-i-1}}_1$ by rotating 
$K^x_0$ about the point $C_{i+1}^{x}$ by an angle $\al$
in the positive direction such that $\sin \al <\sqrt{2h\ep} \cd 2^{-i/2}$,
and the 
rotated copy of $K^x_0$, denoted by $K^{x+2^{-i-1}}_1$ contains the point
$P_{x+2^{-i-1}}$. The same argument shows that 
the distance between the origin and the centre of $K^{x+2^{-i-1}}_1$
is less than $h+3\sqrt{h\ep} \cd \sum_{j=0}^i 2^{-j/2}$.
This defines the circles $K_1^y$ for every 
$y\in D_{i+1} \se D_i$.
It is clear that \eqref{e5} holds for $i+1$ in place of $i$.

Since $\sum_{j=0}^\infty 2^{-j/2}=2+\sqrt 2 <3.2$, 
{\it the distance between the origin and any of the centres of
$K^y_0$ and $K^y_1$ is less than 
\begin{equation}\label{e14}
h+ 9.6 \sqrt{h\ep}< 10\sqrt{h\ep} <\ep
\end{equation} 
for every $y\in D_{i+1}$.}

Now we prove that \eqref{e6} is satisfied with $i+1$ 
in place of $i$.

Let $x\in D_i$, $0<x<1$ be fixed. Let $Q=P_{x-2^{-i-1}}$ and
$K'=K^{x-2^{-i-1}}_0$, $K''=K^{x+2^{-i-1}}_1$. We prove that $K'$ 
and $K''$ intersect each arc
$\iv{A_{j}^0  A_{j}^1}$ $(i<j\le n)$.

The point $C_{i+1}^{x-2^{-i-1}}$ belongs to the arc
$\iv{A_{i+1}^0  A_{i+1}^1}$, and thus it is either inside or on the circle
$K_1$. The point $Q$ is outside $K_1$, and thus the circle
$K'$ must intersect $K_1$ at a point
$S_1$ belonging to the arc $\iv{Q C_{i+1}^{x-2^{-i-1}}}$.  
Let $T_1$ be the other intersection of the circles 
$K'$ and $K_1$. A similar argument shows that the circle
$K'$ intersects $K_0$ at a point
$S_0$ belonging to the arc $\iv{Q C_{i+1}^x}$.  
Let $T_0$ be the other intersection of the circles 
$K'$ and $K_0$. 

We show that $T_0$ and $T_1$ are outside the domain $H$.
Since the distance between the centres of
$K'$ and $K_0$ is less than $h+\ep < 2\ep$, we have 
$$|{S_0 T_0}|\ge 2 \sqrt{1-\ep ^2} >2-\ep> R=\diam H.$$ 
Thus $T_0$ and $S_0$ cannot be both in $H$. 
If $T_0$ is in $H$, then we have $S_0\notin H$, and then
$S_0$ is on the subarc $\iv{P_0 M}$ of $K_0$. Therefore,
$$|{S_0 T_0}|\le |{T_0 M}| + |{MS_0}| \le R+\ep  < 2-\ep  <2\sqrt{1-\ep ^2}$$
which is impossible. Thus $T_0 \notin H$.
The same argument shows $T_1 \notin H$.

If the point $M$ is outside the circle $K'$
then, as $P_1$ is inside $K'$, it follows
that $S_1 \notin H$. If $M$ is inside the circle $K'$, 
then, as $P_0$ is outside $K'$, it follows
that $S_0 \notin H$. If $M$ is on the circle $K'$,
then we have $S_0 =S_1 =M$.
Let $\del H$ denote the boundary of $H$. Then we have 
$$K' \cap (K_0 \cup K_1 )\cap \del H =\{ S\} ,$$
where $S$ is either $S_0$ or $S_1$ (or both if $S_0 =S_1 =M$).

Now $K'$ must
intersect the boundary of $H$ at two points. 
Let the intersection other than $S$ be $U$. Since
$T_0 ,T_1 \notin H$, the point $U$ must be on the arc
$\iv{A_n^0 A_n^1}$. Therefore, the subarc $\iv{C_{i+1}^{x-2^{-i-1}} U}$
of $K'$ is in $H$, and thus 
$K'$ intersects each of the arcs
$\iv{A_j^0 A_j^1}$ for every $i <j\le n$.

A similar argument shows that $K''$ also intersects each of the arcs
$\iv{A_j^0 A_j^1}$ for every $i<j\le n$, and thus \eqref{e6}
is satisfied with $i+1$ in place of $i$.

This completes the construction of the points $C_{i}^x$ and circles
$K^{x}_1$, $K^{x}_0$.

\section{Moving circular arcs} \label{s3}

Let $\si$ be a circular arc of radius $1$ and of arc length less than $1.32$. 
Then the diameter of the set $\si$ is less than $2 \sin(1.32/2) < 1.227=R$.
We describe a continuous motion of $\si$ starting 
from the subarc $\iv{MS_0}$ of $K_0$ lying on the
boundary of $L_1$. Its final position will be the 
subarc $\iv{MS_1}$ of $K_1$ also on the
boundary of $L_1$. 

We choose a large $n$, and take the circles of the sprouting process.
We use the notation of the previous section. The centre of $K_i^x$
is denoted by $O_i^x$ $(i=0,1)$.

First we rotate $\si$ about the point $O_0$
such that its lower endpoint becomes $A_n^0 =C_n^0$. Then we rotate $\si$
about the point $C_n^0$ in the positive direction until it becomes
a subarc of $K_1^{2^{-n}}$. 
Let $\al _n^0$ denote the angle of the rotation. 
Then we rotate $\si$ about the point $O_1^{2^{-n}}$
into $L_0$ with upper endpoint $P_{2^{-n}}$. Then we rotate $\si$
about the point $P_{2^{-n}}$ in the positive direction 
until it becomes
a subarc of $K_0^{2^{-n}}$. 
Let $\be _n^{2^{-n}}$ denote the angle of the rotation. 
Then we rotate $\si$ about the point $O_0^{2^{-n}}$
back into $L_1$ such that its lower endpoint becomes $C_n^{2^{-n}}$. 
Then we rotate $\si$
about the point $C_n^{2^{-n}}$ in the positive direction 
until it becomes
a subarc of $K_1^{2\cd 2^{-n}}$. 
Let $\al _n^{2^{-n}}$ denote the angle of the rotation. 
Then we rotate $\si$ about the point
$O_1^{2\cd 2^{-n}}$
into $L_0$ with upper endpoint $P_{2\cd 2^{-n}}$. Then we rotate $\si$
about the point $P_{2^{-n}}$ in the positive direction
until it becomes
a subarc of $K_0^{2\cd 2^{-n}}$. 
Let $\be _n^{2\cd 2^{-n}}$ denote the angle of the rotation. 
Then we rotate $\si$ about the point $O_0^{2\cd 2^{-n}}$
back into $L_1$ such that its lower endpoint becomes $C_n^{2\cd 2^{-n}}$. 
We continue this process until $\si$ becomes a subarc of
$K_0^{(2^n -1)\cd 2^{-n}}$ with lower endpoint $C_n^{(2^n -1)\cd 2^{-n}} =A_n^1$.
Then we rotate $\si$
about the point $A_n^1$ in the positive direction
until it becomes a subarc of $K_1$. 
Let $\al _n^{(2^n -1)\cd 2^{-n}}$ denote the angle of the rotation. 
Then we rotate $\si$ about the point $O_1$
until it reaches its final position $\iv{MS_1}$.

Let $x\in D_n$, $x<1$.
Let $H_i^x$ denote the horn shaped domain bounded by the
subarc $\iv{P_x C_i^x}$ of $K_0^x$, the
subarc $\iv{P_{x+2^{-i}} C_i^{x}}$ of $K_1^{x+2^{-i}}$, and
the arc $\iv{P_x P_{x+2^{-i}}}$. We put $\De (h) =H_0^0$. Then
$\De (h)$ is the horn shaped domain bounded by the
subarc $\iv{P_0 M}$ of $K_0$, $\iv{P_1 M}$ of $K_1$ and the
arc $\iv{P_0 P_1}$. We put
$$T_n (h)= \bigcup_{x\in D_n , \ x<1} H_n^x .$$
Our next aim is to estimate the area of $T_n (h)\se \De (h)$.
The area (Lebesgue measure) of a set $A\su \sik$ will be
denoted by $m(A)$.

If a point $P$ belongs to $T_n (h) \se \De (h)$, then
there is a smallest $j>0$ such that $P\in H_j^x$ for some $x\in D_j$,
$x<1$. If $j=i+1$, then 
$$P\in H_{i+1}^x \se \bigcup_{y\in D_i , \ y<1} H_i^y$$
for some $x\in D_{i+1}$, $x<1$. Now we have either $x\in D_i$
or $x-2^{-i-1} \in D_i$. In the first case $P\in H_{i+1}^x \se H_i^x$,
while in the second case $P\in H_{i+1}^x \se H_i^{x-2^{-i-1}}$.
We have proved the following:
\begin{equation}\label{e20}
\begin{split}
T_n (h)\se \De (h) \su & \bigcup_{i=0}^{n-1} \bigcup_{x\in D_i \ x<1}
\left( H_{i+1}^x \se H_i^x \right) \cup \\
& \cup \bigcup_{i=0}^{n-1} \bigcup_{x\in D_{i+1} \se D_i}
\left( H_{i+1}^x \se H_i^{x-2^{-i-1}} \right) .
\end{split}
\end{equation} 
\begin{figure}
\centering
\includegraphics[width=4in]{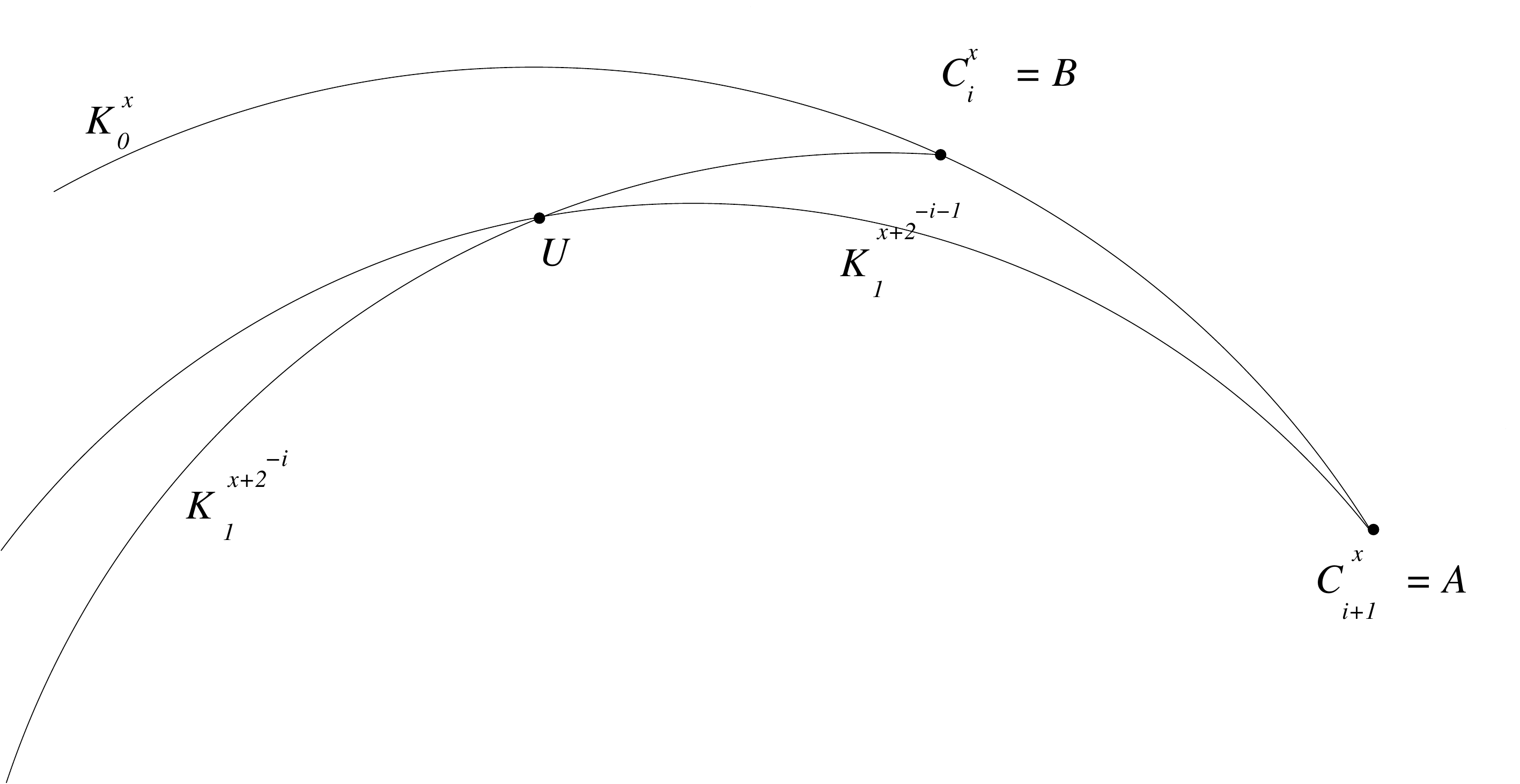}
\caption{$H_{i+1}^x \se H_i^x$}
\label{fig13}
\end{figure}

The set $H_{i+1}^x \se H_i^x$ is bounded 
by three subarcs of the circles $K_0^x , K_1^{x+2^{-i}}$ and 
$K_1^{x+2^{-i-1}}$ (see Figure \ref{fig13}). Let $U$ denote the intersection 
of $K_1^{x+2^{-i}}$ and $K_1^{x+2^{-i-1}}$ belonging to $L_1$.
We want to apply Lemma \ref{l3} with the choice 
$K=K_0^x$, $A=C_{i+1}^x$, $B=C_i^x$,
$\al =\al_{i+1}^x$ being the angle of rotation about $A$ mapping 
$K_1^{x+2^{-i-1}}$ onto $K_0^x$, and $\be =\al_{i}^x$ being the angle of rotation 
about $B$ mapping 
$K_1^{x+2^{-i}}$ onto $K_0^x$. Now we check that (for $n$ large enough)
the conditions
of Lemma \ref{l3} are satisfied. First we need the following 
lemma (its proof will be given in Section \ref{s5}).
\begin{lemma}\label{l8}
There exists a positive constant $c$ depending on $h$ and $\ep$
but not on $n$ such that $|C_i^x C_{i+1}^x |\le c/n$
and $|C_i^x C_{i+1}^{x+2^{-i-1}} |\le c/n$ hold
for every $n=1,2,\ldots$, $0\le i<n$ and $x\in D_i$, $x<1$.
\end{lemma}
By this lemma, we have 
\begin{equation}\label{e21}
\eta =C_{i+1}^x O_0^x C_i^x \angle <2c/n<1/5
\end{equation} 
if $n>10c$.
Next we need the following estimates for $\al _i^x$.
\begin{lemma}\label{l10}
For every $1\le i\le n$ and $x\in D_i$, $x<1$ we have 
\begin{equation}\label{e22}
0.9 \cd \frac{h\ep \cd 2^{-i}}{t\cd \sqrt{1-(t/2)^2}} <\al _i^x <
1.1 \cd \frac{h\ep \cd 2^{-i}}{t\cd \sqrt{1-(t/2)^2}} ,
\end{equation} 
where $t=|P_{x+2^{-i}} C_i^x |$.
\end{lemma}
We prove Lemma \ref{l10} in Section \ref{s5}. 

Clearly, $|C_i^x C_{i+1}^{x} |\ge 1/n$, and thus
$\eta =C_{i+1}^x O_0^x C_i^x \angle > 1/n.$
Since $t=t_i >\ep$ and $t\le \ep +R \leq 2-4\ep$, \eqref{e22} gives
\begin{equation}\label{e34}
\al _i^x \le 1.1 \cd \frac{h\ep \cd 2^{-i}}{\ep \cd \sqrt{1-(1-2\ep)^2}} \le
1.1 \cd \frac{h \cd 2^{-i}}{\sqrt{2\ep}} <2^{-i}\cd 10^{-3} < 2^{-i},
\end{equation} 
where we used that $\ep < 1/2$ and $h/\ep<10^{-3}$.
If $i> \log _2 n$, then this gives $\al _i^x  <1/n<\eta$.
Thus, when applying 
Lemma \ref{l3}, the condition $\be <\eta$ is satisfied if 
$i> \log _2 n$. Finally, we prove that the condition $\al /\be <3/4$
is satisfied if $i$ is large enough. We have, by \eqref{e22},
\begin{equation}\label{e32}
\al /\be =\al_{i+1}^x /\al_{i}^x <\ha \cd \frac{1.1}{0.9}\cd
\frac{t_i \cd \sqrt{1-(t_i /2)^2}}{t_{i+1}\cd \sqrt{1-(t_{i+1}/2)^2}} ,
\end{equation} 
where $t_i =|P_{x+2^{-i}} C_i^x |$ and $t_{i+1} =|P_{x+2^{-i-1}} C_{i+1}^x |$.
Clearly, 
\begin{equation}\label{e33}
|t_{i+1}-t_i |\le  |P_{x+2^{-i}} P_{x+2^{-i-1}}|+|C_i^x C_{i+1}^x |<
2^{-i}+(c/n).
\end{equation} 
Since the function $t\mapsto t\cd \sqrt{1-(t/2)^2}$ is positive and
uniformly continuous on the interval $[0,2-4\ep ]$, there exists
a positive number $v$ only depending on $\ep$ such that 
$$(t\cd \sqrt{1-(t/2)^2})/(t'\cd \sqrt{1-(t'/2)^2} ) < 1.1$$
whenever $t,t'\in [0,2-4\ep ]$ and $|t-t' |<v$.
If $n> (1+c)/v$, and $i> \log _2 n $, then \eqref{e32} and \eqref{e33}
imply $\al /\be <(1/2)\cd (1.1/0.9)\cd 1.1 <3/4$.

Therefore, assuming Lemmas \ref{l8} and \ref{l10} we can see 
that all conditions of Lemma \ref{l3} are satisfied 
provided that $n>\max (10c, (1+c)/v)$ and $i>\log _2 n $. 
Therefore, the circles $K_1^{x+2^{-i}}$ and $K_1^{x+2^{-i-1}}$ intersect at a point
$U$ such that $|UC_{i+1}^x|<20\eta <40c/n$. 
Rotating the subarc $\iv{UC_{i+1}^x}$ of $K_1^{x+2^{-i-1}}$ about the point
$C_{i+1}^x$ in the negative direction by angle $\al _{i+1}^x$ we obtain 
a horn shaped domain $H'$ which contains $H_{i+1}^x \se H_i^x$. 
Then 
$$m(H')=|UC_{i+1}^x|^2 \cd \al _{i+1}^x /2 < 
(40c /n)^2 \cd 2^{-i}=1600c^2 \cd 2^{-i}/n^2$$
by \eqref{e34} and Lemma \ref{l6}. We find that 
$$m(H_{i+1}^x \se H_i^x ) < 1600c^2 \cd 2^{-i}/n^2$$ 
if $n>\max (10c, (1+c)/v)$ and $i>\log _2 n$. 
A similar argument shows that if if $n>\max (10c, (1+c)/v)$ and $i>\log _2 n$, 
then
$$m(H_{i+1}^x \se H_i^{x-2^{-i-1}})< 10000c^2 \cd 2^{-i}/n^2 .$$ 
The only difference is that we have to apply Lemma \ref{l9}
instead of Lemma \ref{l3}. We omit the details.
Since the cardinality of $D_i$ equals $2^i +1$, we obtain that the area of the set
\begin{equation*}
\bigcup_{i= \lceil \log _2 n \rceil}^{n-1} \bigcup_{x\in D_i \ x<1}
\left( H_{i+1}^x \se H_i^x \right)
 \cup \bigcup_{i=\lceil \log _2 n \rceil}^{n-1} \bigcup_{x\in D_{i+1} \se D_i}
\left( H_{i+1}^x \se H_i^{x-2^{-i-1}} \right) 
\end{equation*} 
is less than
$$2\cd n\cd 2^i \cd 10000c^2 \cd 2^{-i}/n^2 =20000c^2/n$$
if $n>\max (10c, (1+c)/v)$. If $i\le \log _2 n$, 
then for all $x \in D_i, x<1$, $H_i^x$ is contained in the disc
of centre $M$ and radius $iR/n <2(\log _2 n)/n$. Thus we have that the set

\begin{equation*}
\bigcup_{i=0 }^{\lfloor \log _2 n \rfloor} \bigcup_{x\in D_i \ x<1}
\left( H_{i+1}^x \se H_i^x \right)
\cup \bigcup_{i=0}^{\lfloor \log _2 n \rfloor} \bigcup_{x\in D_{i+1} \se D_i}
\left( H_{i+1}^x \se H_i^{x-2^{-i-1}} \right) 
\end{equation*} 
is contained in the disc
of centre $M$ and radius less than $2((\log _2 n) +1)/n$.  
Thus the set $T_n (h) \se \De (h)$ can be covered by a disc of radius
$2((\log _2 n) +1)/n$ and a set of area $20000c^2 /n$. Since $c$ does not
depend on $n$ and
$(\log _2 n)/n \to 0$ as $n\to \infty$, we have proved the following.
\begin{lemma}\label{l11}
The area of the set $T_n (h) \se \De (h)$ tends to zero as $n\to \infty$.
\end{lemma}

\section{Proof of Theorem \ref{t2}} \label{s4}

Let $\si _0$ be a circular arc of radius $1$ and of arc length less than $1.32$, 
and let $\si_1$ be obtained 
from $\si _0$ by rotating it about one of its end points by a given angle
$h \in (0,\pi )$. Let $\ep >0$ be given.
We have to prove that there is a continuous motion
bringing $\si _0$ to $\si _1$ 
such that the area touched by the
moving arc is less than $\ep$. Let $\diam \si _0$ denote 
the diameter of $\si _0$, then we have 
$\diam \si _0 < 1.227$. We may assume that $\ep <10^{-6}$.
In the argument below we fix $\ep$ with this property.

We shall use the notation of the previous two sections. 
We shall indicate the dependence of the objects
on the angle $h$ whenever necessary. 
That is, $K_0 (h)$ and $K_1 (h)$ will denote
the circles of radius $1$ and centres $O_0 =(-\sin (h/2),0)$ and 
$O_1 =(\sin (h/2),0)$, respectively, for every $h\in (0,\pi )$. 
The point with coordinates $(0,\cos (h/2))$
is denoted by $M(h)$. 

Let $\ol K _i (h)$ denote the disc bounded by $K_i (h)$ $(i=0,1)$.
The closures of the sets $\ol K _0 (h) \se \ol K _1 (h)$ and $\ol K _1 (h)\se 
\ol K _0 (h)$ are denoted by $L_0 (h)$ and $L_1 (h)$. 

We may assume that $\si _0 =\iv{M(h) S_0 (h)}$ and 
$\si _1 =\iv{M(h) S_1 (h)}$, where $\iv{M(h) S_0 (h)}$ is 
a subarc of $K_0 (h) \cap L_1 (h)$ and $\iv{M(h) S_1 (h)}$ is
a subarc of $K_1 (h)$.

The circle of radius $\ep$ and centre $M(h)$ intersects
the lune $L_0 (h)$ in the arc $\iv{P_0 (h)  P_1 (h)}$, 
where $P_0 (h)\in K_0 (h)$ and $P_1 (h)\in K_1 (h)$.

The intersection of $L_0 (h)$ and the disc of centre $M(h)$ and radius $\ep$
is denoted by $\De (h)$. Clearly, $m(\De (h))< \ep ^2 \pi <\ep /2$.

Let $\La$ denote the set of all positive numbers $\la$ having the following
property: for every $0<h<\pi$ there is a 
set $\Si (h, \la )$ such that
$$m(\Si (h,\la ) \se \De (h))< \la \cd h,$$
and an arc can be moved 
continuously in the set $\Si (h, \la )$ with initial
position $\iv{M(h) S_0 (h)}$ and final position $\iv{M(h) S_1 (h)}$.

If we rotate $\iv{M(h) S_0 (h)}$ about the point $M(h)$ in the 
positive direction by angle $h$, then the area of the set $\Si$
of points touched by the moving arc equals $(\diam \si _0 )^2 \cd h/2$ 
by Lemma \ref{l6}. Therefore, we have 
$(\diam \si _0 )^2/2 \in \La$, and thus $\La \ne \emp$.

Our aim is to show that if $\la \in \La$, then $\frac{q+3}{4}\cd \la \in \La$ 
(see Lemma \ref{l13} below), where $0<q<1$ is a constant depending only on $\ep$. 

We will need the following lemma about the sum of the 
angles $\beta_n^{x}$ of rotations around the points $P_x$.

\begin{lemma}\label{l+}
There exists a positive constant $q < 1$ depending only on $\ep$ such that 
\begin{equation}\label{eqbeta}
\sum_{x \in D_n, \ 0<x<1} \beta_n^{x} \leq q \cd h.
\end{equation}

\end{lemma}

In Section \ref{s5} we prove that the statement of Lemma \ref{l+} 
is true with $q=1-\ep^4 $.

We shall use the notation
$$U(A,\eta )=\{ x\in \sik : \dist (x,A)<\eta \}$$
for every $A\su \sik$ and $\eta >0$.
\begin{lemma}\label{l12}
For every $0<h<\ep ^2 /100$ 
and $\la \in \La$ there is a set $\Si ' (h,\la )$
such that
\begin{equation}\label{e25}
m\left( \Si ' (h,\la )\se U(\De (h), \sqrt{h} ) \right)
< \frac{1+q}{2} \cd \la \cd h ,
\end{equation} 
and an arc can be moved 
continuously in the set $\Si '(h, \la )$ with initial
position $\iv{M(h) S_0 (h)}$ and final position $\iv{M(h) S_1 (h)}$.
\end{lemma}

\proof
Let $0<h<\ep ^2 /100$ and $\la \in \La$ 
be fixed. In the argument below we shall
suppress the reference to $h$; that is, we shall write $K_0$
instead of $K_0 (h)$ etc.

We apply the sprouting process
as described in Section \ref{s2}, and construct the set
$T_n$ as described in Section \ref{s3}. By Lemma \ref{l11},
we can choose an $n$ such that 
\begin{equation}\label{e36}
m(T_n \se \De (h))< \frac{1-q}{2} \cd \la \cd h
\end{equation} 
holds. The movement of the arc consists of rotating 
about the points $C_n ^x$ and $P_x$ alternately, and 
moving the arc along some circles
in between. The angle of the rotation about the point $P_x$
equals $\be _n^x$. The aim of this rotation is
to bring the arc from $K_1^x$ to $K_0^x$.

Now we shall replace the set touched by the arc during the rotation
about the point $P_x$ by a smaller set using the condition
$\la \in \La$. We fix $x\in D_n$, $0 <x <1$, and write $\be$ in place of $\be _n^x$.
By definition we can move an arc continuously in 
the set $\Si (\be ,\la )$ such that its initial
position is $\iv{M(\be ) S_0 (\be )}$ and its final position
is $\iv{M(\be ) S_1 (\be )}$. 

Let $\tau _x$ denote the isometry
mapping $O_0 (\be )$ into $O_0^x$ and $O_1 (\be )$ into $O_1 ^x$.
Note that $|O_0 (\be )O_1 (\be )|=|O_0^x O_1^x |=2\sin (\be /2)$.
The isometry $\tau _x$ is obtained by a reflection about the $y$
axis followed by a suitable rotation by an angle close to $\ep$
which maps the point $M(\be )$ into $P_x$.
Then $\tau _x$ is an isometry mapping the circle $K_1 (\be )$ into $K_1 ^x$
and the circle $K_0 (\be )$ into $K_0 ^x$. 
It is clear that we can move an arc continuously in 
the set $\tau _x \left( \Si (\be ,\la )\right)$ such that its initial
position is a subarc of $K_1^x \cap L_0 (h)$ and its final position
is a subarc of $K_0^x \cap L_0 (h)$.

\begin{figure}
\centering
\includegraphics[width=4in]{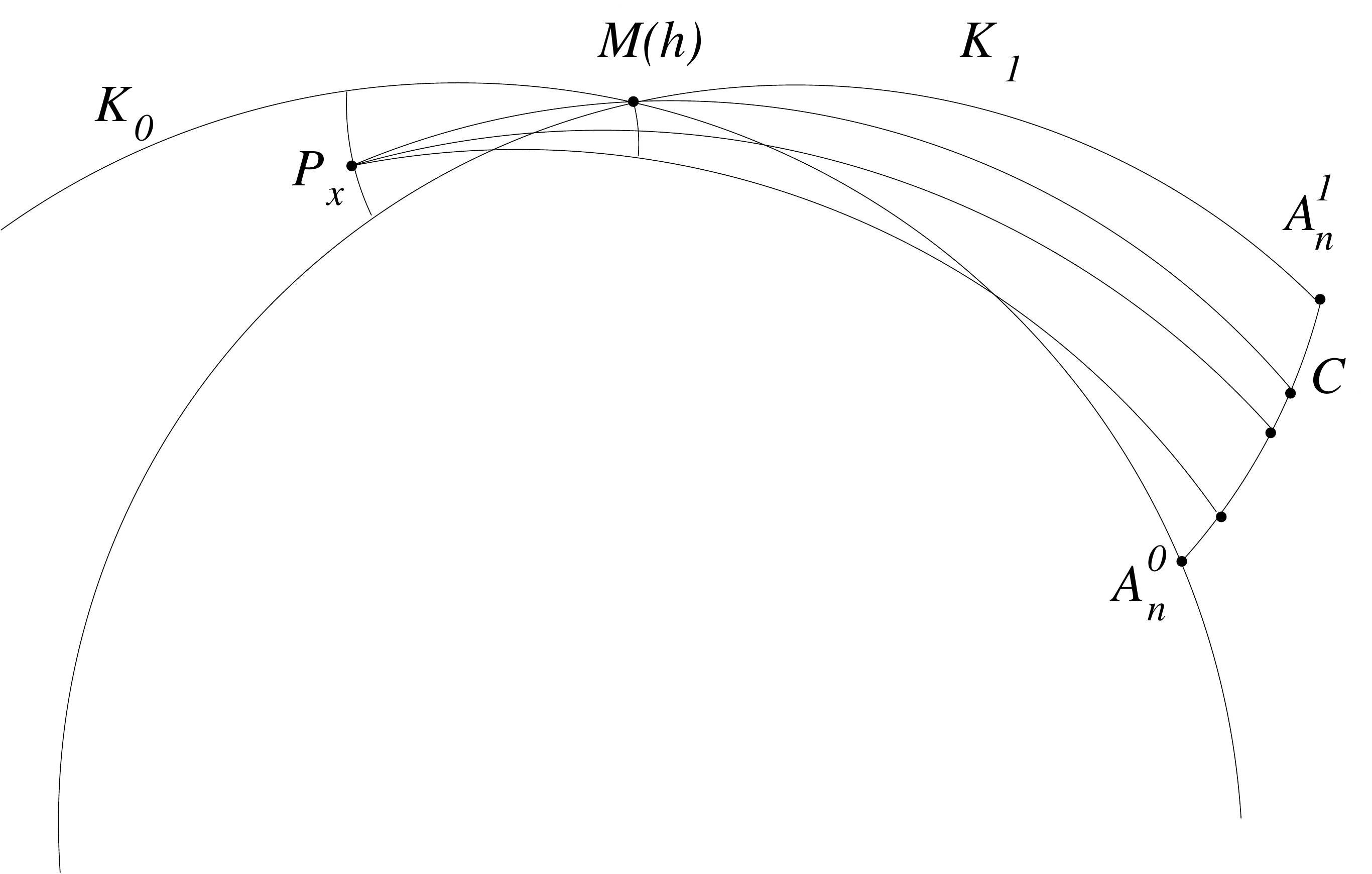}
\caption{Proof of Lemma \ref{l12}}
\label{fig10}
\end{figure}

Let $\Si ' (h,\la )$ denote the set 
\begin{equation}\label{e27}
T_n  \cup \bigcup_{x\in D_n , \ 0 <x<1}
\tau _x \left( \Si (\be _n^x ,\la )\right) .
\end{equation} 
We prove that $\Si ' (h,\la )$ satisfies the condition of the lemma.
It follows from the construction that we can move an arc 
continuously in $\Si ' (h,\la )$ from the initial
position $\iv{M S_0}$ to the final position $\iv{M S_1}$.
In order to prove \eqref{e25} first we show that 
\begin{equation}\label{e26}
\tau _x (\De (\be _n^x )) \su U(\De (h), \sqrt h )
\end{equation} 
for every $x\in D_n$, $0 <x<1$. The set 
$\tau _x (\De (\be _n^x ))$ is bounded by subarcs of $K_1^x$, $K_0^x$
and of the circle with centre $P_x$ and radius $\ep$ (see Figure
\ref{fig10}).
Let $K'_0$ denote the circle of radius $1$ going through
the points $P_x$ and $M=M(h)$ obtained from $K_0$ by rotating
it about the point $M$ in the positive direction by an angle $<h$.
It is clear that $K'_0$ intersects the arc $\iv{A_n^0 A_n^1}$.
Let $C$ denote the point of intersection. The circle $K_0^x$
also intersects the arc $\iv{A_n^0 A_n^1}$ at the point $C_n^x$.
Since $|C_n^x C|<|A_n^0 A_n^1 |<R\cd h<2h$ and $|P_x C|>\ep >3\sqrt{2h}$,
$|P_x C|<R+\ep =2-4\ep <2-2h$, we may apply Lemma $\ref{l1}$
and find that $K_0^x$ is obtained from $K'_0$ by rotating it about
the point $P_x$ by an angle $\al$ such that $\sin \al <2\sqrt{2h}$.
Thus $\al <4\sqrt{2h}$.
The same argument shows that 
$K_1^x$ is obtained from $K'_0$ by rotating it about
the point $P_x$ by an angle $<4\sqrt{2h}$.

Since $|P_x M|=\ep$, we can see that 
every point $A$ of $\tau _x (\De (\be _n^x ))$ can be obtained from
a point $B$ of the subarc $\iv{P_x M}$ of $K'_0$
by rotating it about the point $P_x$ by an angle $<4\sqrt{2h}$.
Then $|AB|<\ep \cd 4\sqrt{2h}<6\ep \sqrt h < \sqrt h$. Since $B\in \De (h)$,
this proves \eqref{e26}.

We have
$$m\left( \Si (\be _n ^x ,\la ) \se \De (\be _n ^x ) \right)
<\la \cd \be _ n ^x$$ 
for every $x\in D_n$, $0 <x<1$ by the definition of $\la$. 
Then it follows from \eqref{e36}, \eqref{e27} and \eqref{e26} that 
the area of the set $\Si '(h,\la ) \se U(\De (h), \sqrt h )$ is less than
$$\left(\frac{1-q}{2} \cd \la \cd h \right)      +\sum_{x\in D_n ,\ 0 <x<1} \la \cd \be _n^x .$$
Since $\sum_{x\in D_n ,\ x>0} \be _n^x \le q \cd h$ by Lemma \ref{l+},
we find that the area in question is less than
$$\frac{1-q}{2} \cd \la \cd h+q \cd \la \cd h =\frac{1+q}{2} \cd \la \cd h.$$
This completes the proof of Lemma \ref{l12}.
\hfill $\square$

\begin{lemma}\label{l13}
For every $\la \in \La$ we have $\frac{q+3}{4} \cd \la \in \La$.
\end{lemma}

\proof
Let $h\in (0,\pi )$ and $\la \in \La$ be arbitrary. 
Let $G$ be an open set containing $\De (h)$ such that
$m(G\se \De (h))< \frac{1-q}{4} \cd \la \cd h$.
Since $\De (h)$ is compact, there is an $\eta >0$ such that
$U(\De (h), \eta )\su G$. 
We fix a positive number $\de$ such that $\de <\min (\ep ^2 /100 , \eta ^2 )$. 

\begin{figure}
\centering
\includegraphics[width=4in]{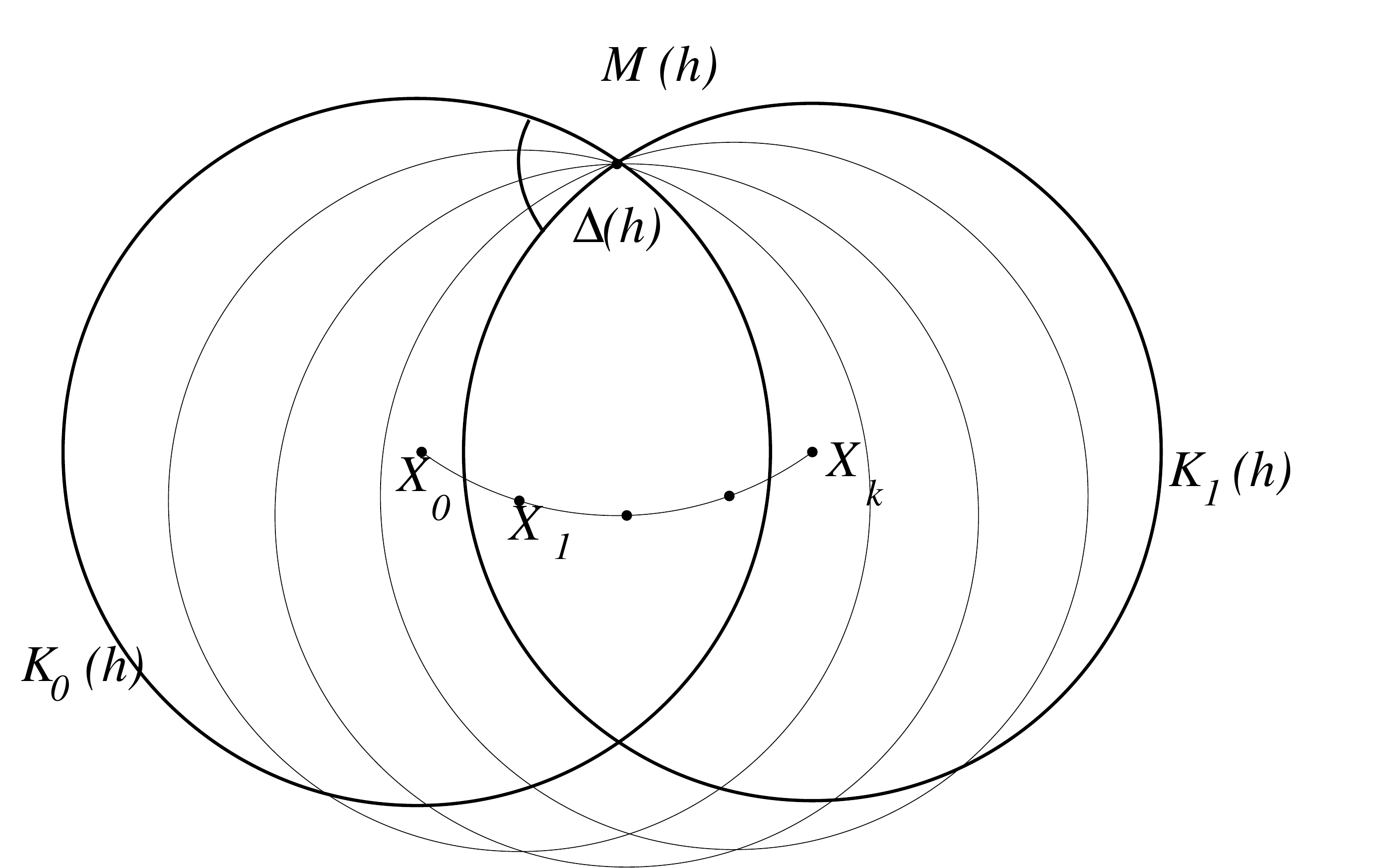}
\caption{Proof of Lemma \ref{l13}}
\label{fig9}
\end{figure}

We consider the subarc $\iv{O_0 O_1}$ of the circle of centre $M(h)$
and radius $1$. We select points of division $O_0 =X_0 ,X_1 \stb X_k =O_1$
on the arc $\iv{O_0 O_1}$ such that each of the angles
$h_i =X_i M(h) X_{i-1} \angle$
$(i=1\stb k)$ is less than $\de$ (see Figure \ref{fig9}).
Note that $h_1 +\ldots +h_k =h$.

The circle of centre $X_i$ and radius $1$ will be denoted by $N_i$
$(i=0\stb m)$. Thus $N_0 =K_0 (h)$ and $N_k =K_1 (h)$.

For every $i=1\stb k$ there is an isometry (in fact, a rotation
by an angle less than $h$)
$\kappa _i$ mapping $K_0 (h_i )$ onto $N_{i-1}$ and $K_1 (h_i )$
onto $N_i$. It is easy to see that $\kappa _i (\De (h_i )) \su \De (h)$
for every $i=1\stb k$. We put
$$\Si =\sumik \kappa _i \left( \Si ' (h_i ,\la ) \right) ,$$
where the set $\Si ' (h_i ,\la )$ is given by Lemma \ref{l12}.
It is clear that we can move an arc 
continuously in the set $\Si$ with initial
position $\iv{M(h) S_0 (h)}$ and final position $\iv{M(h) S_1 (h)}$.
Since $\sqrt{h_i} <\sqrt \de <\eta$ for every $i$, it follows from
\eqref{e25} that
$$m\left( \kappa _i \left( \Si ' (h_i ,\la ) \right)
\se U(\De (h), \eta )\right) < \frac{1+q}{2} \cd \la \cd h_i$$
for every $i$. Therefore, we have
$$m\left( \Si \se U(\De (h), \eta ) \right) < 
\sumik \frac{1+q}{2} \cd \la \cd h_i =\frac{1+q}{2} \cd \la \cd h.$$
Hence, by $m(U(\De (h), \eta ) \se \De (h))<\frac{1-q}{4}  \cd \la \cd h$
we obtain 
$$m\left( \Si \se \De (h)\right) < 
\left(\frac{1-q}{4} + \frac{1+q}{2} \right) \cd \la \cd h= \frac{q+3}{4} \cd \la \cd h.$$
Since $0<h<\pi$ was arbitrary, this means,  by the definition of $\La$
that $\frac{q+3}{4} \cd \la \in \La$, which completes the proof of the lemma.
\hfill $\square$

\bsk
In possession of Lemma \ref{l13} we prove Theorem \ref{t2}
as follows. 
Let $\si _0$ and $\si _1$ be as in Theorem \ref{t2}. We may assume that
the angle $h$ of the rotation bringing $\si _0$ onto $\si _1$ is less than
$\pi$. Then we may also assume that $\si _0 =\iv{M(h)S_0 (h)}$ and
$\si _1 =\iv{M(h)S_1 (h)}$.

As $\La \ne \emp$, we can choose a positive number $\la _0 \in \La$.
Then, by Lemma \ref{l13}, $\la _n =\left(\frac{q+3}{4}\right)^n \cd \la _0 \in \La$
for every $n$. Since $0<q<1$, thus $0<\frac{q+3}{4} <1$ and $\la _n \to 0$ as $n\to \infty$, 
hence we can choose a $\la \in \La$ such that $\la  <\ep /(2\pi )$. 

By the definition of $\La$
this means that there is a set $\Si (h, \la )$ such that
$m(\Si (h,\la ) \se \De (h))< \la  \cd h <\ep /2$,
and an arc can be moved 
continuously in the set $\Si (h, \la )$ with initial
position $\iv{M(h) S_0 (h)}$ and final position $\iv{M(h) S_1 (h)}$.
Since $m(\De (h))<\ep /2$, we have $m(\Si (h, \la ))<\ep$,
and thus Theorem \ref{t2} is proved, 
subject to the Lemmas
whose proofs were postponed. 
\hfill $\square$

\section{Proof of the lemmas} \label{s5}

{\bf Proof of Lemma \ref{l1}.} 
We shall use the following notation. 
Let $A$ be the closest point of $K$ to $P$, and 
let $B$ be the point obtained by reflecting $A$ about $O$.
Then the points $B,O,P,A$ are collinear, and $|PA|=d$.
Let $R$ be the point of the segment $OA$ such that $|RA|=2d$.
Let $C,D \in K$ be such that the segment $CD$ is perpendicular to $OA$ and
contains $R$. We assume that $C$ and $Q$ belong to the same semicircle
with endpoints $B$ and $A$.
Then we have
\begin{equation}\label{e1}
|RC|=\sqrt{|RA|\cd |RB|}= \sqrt{2d\cd (2-2d)} =2\sqrt{d-d^2} .
\end{equation}
Now we turn to the proof of the lemma.

\begin{figure}
\centering
\includegraphics[width=3in]{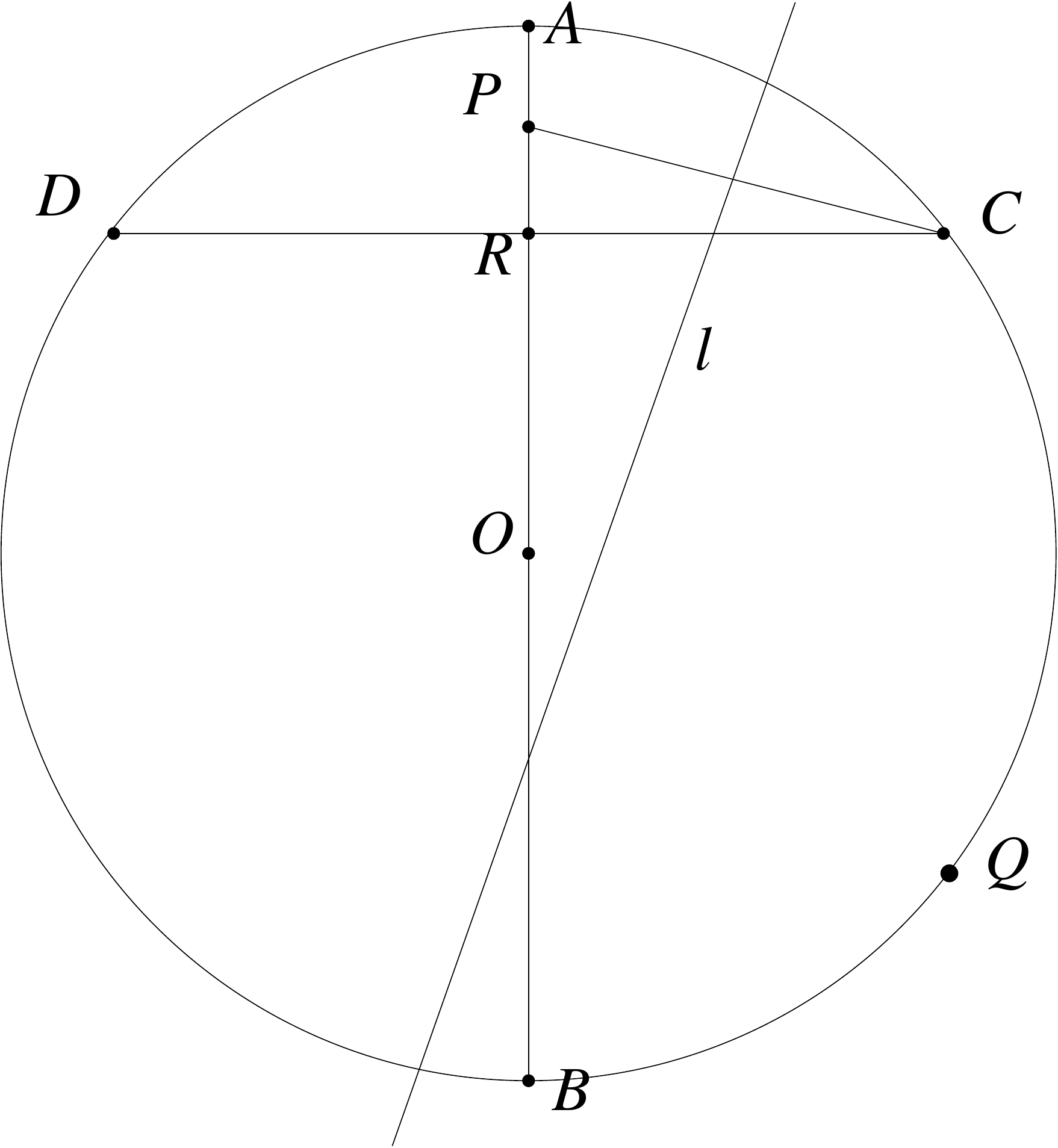}
\caption{Case (i)}
\label{fig:in}
\end{figure}

\bsk \noi
(i) Suppose that $P$ is in the interior of $K$. 
Then we have 
$$|PC| =\sqrt{|PR|^2 +|RC|^2} =\sqrt{d^2 + 4\cd (d-d^2 )}<2\sqrt d.$$ 
Since $|QP|\ge 2\sqrt d$, it follows 
that $Q$ belongs to the arc $\iv{CB}$.
By \eqref{e1} we have
$$|BC|^2 =|BR|^2 +|RC|^2=(2-2d)^2 + 4(d-d^2 )=4-4d<(2-d)^2 =|BP|^2 ,$$
and thus $|BP|>|BC|$. Now we prove $|QP|>|QC|$.

Let $\ell$ denote the 
perpendicular bisector of the segment $PC$. Since
$|BP|>|BC|$ and $|DP|<|DC|$, 
it follows that $\ell$ intersects the arc $\iv{BD}$.
It is clear that
$\ell$ intersects the arc $\iv{AC}$ as well, and thus the arc $\iv{BC}$
is disjoint from $\ell$. Then every point of 
$\iv{BC}$ is closer to $C$ then to $P$. In particular, $|QP| >|QC|$. 
Thus 
the circle of radius $|QP|$ and centre $Q$ intersects 
the arc $\iv{AC}$ at a unique point $P'$ such that 
the orientation of the triangle $\De_{PQP'}$ is positive. Let $\al =P'QP\angle$.
It is clear that $\al <\pi /2$, and 
rotating $K$ about the point $Q$ in the positive direction by
angle $\al$, the circle obtained will contain $P$. We show
$\sin \al <2\sqrt{d}$. Indeed, we have
$$\al <CQP\angle <CQD\angle =CBD\angle =\tfrac{1}{2} \cd COD\angle =
COR\angle ,$$
and thus $\sin \al < \sin COR\angle <2\sqrt{d}$.

\msk \noi
(ii) Suppose that $P$ is outside the circle $K$. 
Then we have 
$$|PC|=\sqrt{(3d)^2 +(4d-4d^2 )}<\sqrt{9d}=3\sqrt d .$$
Since $|QP|\ge 3\sqrt{d}$, it follows that $Q$ belongs to the arc $\iv{CB}$.
Let $E$ be the point obtained by reflecting $D$ about $O$.
Then $\De_{DEC}$ is a right triangle, and $CE$ is parallel to $AB$.
Let $S$ be the point obtained by reflecting $R$ about $O$. Then
we have 
$$|PS|=|PR|+2\cd |RO|=3d+2(1-2d)=2-d,$$ 
and thus
$$|PE|^2 =|PS|^2 +|SE|^2 =(2-d)^2 +4\cd (d-d^2 )= 4-3d^2 .$$
This gives
$$2-d < |PE|<2=|DE|.$$
\begin{figure}
\centering
\includegraphics[width=3in]{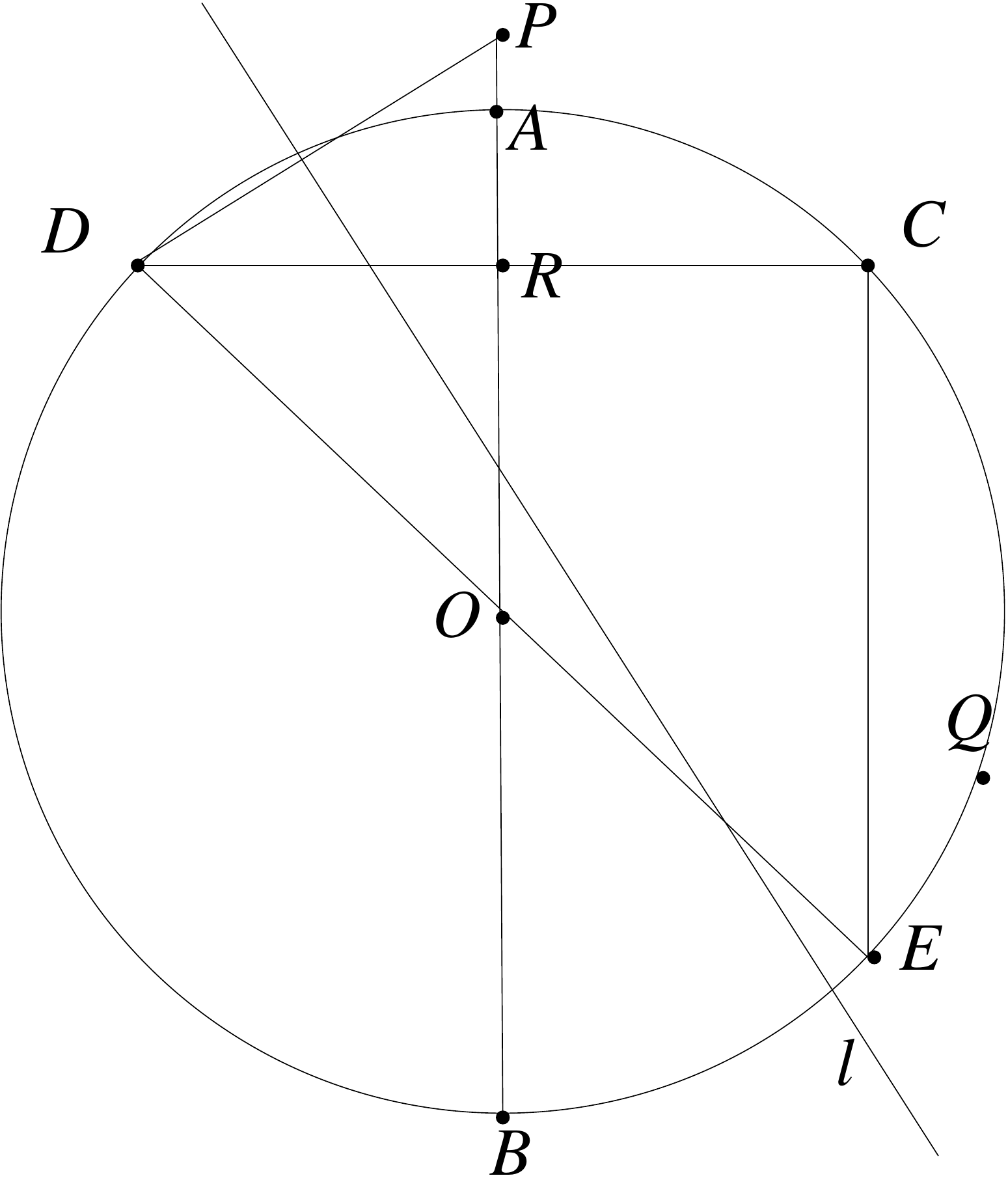}
\caption{Case (ii)}
\label{fig:out}
\end{figure}
Then, by $|PQ|\le 2-d$ it follows that $Q$ belongs to the arc $\iv{CE}$.
We prove that $|QP|< |QD|$. Let $\ell$ denote the 
perpendicular bisector of the segment $DP$. Since
$|PE|<|DE|$ and $|PB|>|DB|$, 
it follows that $\ell$ intersects the arc $\iv{BE}$.
Similarly, $|AP|=d<3d<|DP|$ implies that  
$\ell$ intersects the arc $\iv{DA}$ as well. Therefore, the arc $\iv{CE}$
is disjoint from $\ell$, and thus every point of 
$\iv{CE}$ is closer to $P$ then to $D$. This proves
$|QP|< |QD|$. 
Therefore,
the circle of radius $|QP|$ and centre $Q$ intersects 
the arc $\iv{DA}$ at a unique point $P'$ such that 
the orientation of the triangle $\De_{PQP'}$ is negative. Let $\al =PQP'\angle$.
It is clear that $\al <\pi /2$, and
rotating $K$ about the point $Q$ in the negative direction by
angle $\al$, the circle obtained will contain $P$. 

It is easy to check that the tangent line of $K$ at the point $C$
does not separate the points $P$ and $Q$. Therefore, we have
$PQP'\angle \le CQD \angle$, and thus
$$\al =PQP'\angle \le CQD \angle =CBD\angle =\tfrac{1}{2} \cd COD \angle = COR \angle .$$
This gives $\sin \al <\sin COR \angle <2\sqrt d$.

We have proved that $\sin \al <2\sqrt d$ in both cases.
Let $O'$ denote the centre of $K'$. Then $\De_{OQO'}$ is an isosceles
triangle such that $|OP|=|O'P|=1$ and $OPO' \angle =\al$.
Thus $|OO'|=2\sin (\al /2)<2\sin \al <4\sqrt d $.
\hfill $\square$

\bsk \noi
{\bf Proof of Lemma \ref{l3}.}
\begin{figure}
\centering
\includegraphics[width=3in]{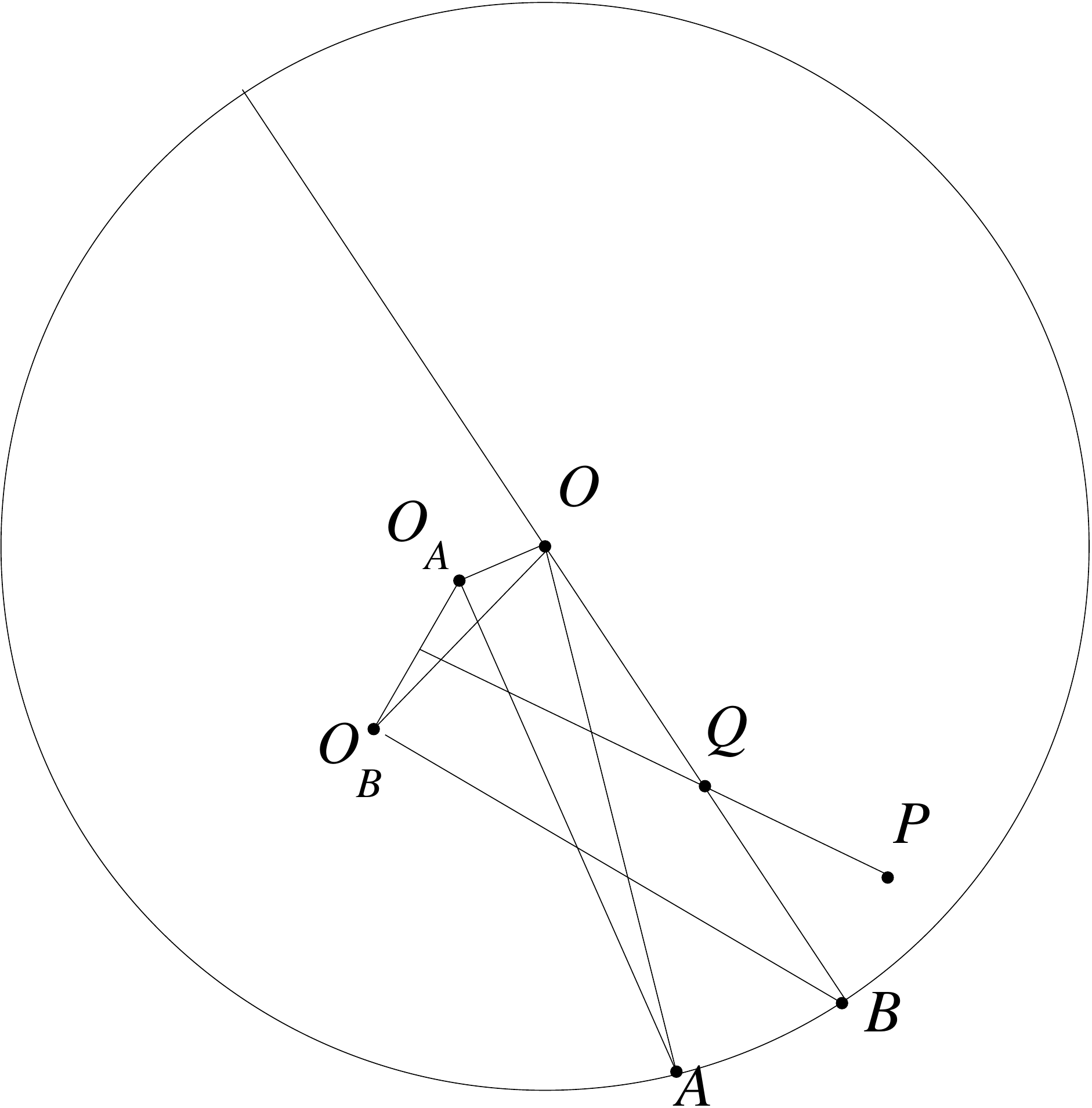}
\caption{Lemma \ref{l3}}
\label{fig:inout}
\end{figure}
Let $O_A$ and $O_B $ denote the centres of $K_A$ and $K_B$, respectively.
Then we have $|{O_AO}|=2\sin (\al /2) < 2\sin (\be /2) =|{O_B O}|$, and
\begin{align*}
O_AOO_B \angle & =O_A O A \angle +AOB \angle -BOO_B \angle = \\
& =\left(  \frac{\pi}{2}-\frac{\al}{2} \right) +\eta - 
\left( \frac{\pi}{2}- \frac{\be}{2} \right) =\eta +\frac{\be -\al}{2} .
\end{align*}
Since $\eta +(\be -\al )/2 >0$, we obtain that the orientation of
the triangle $OO_A O_B$ is positive.
If $OO_AO_B \angle =\ga$, then 
$$O_AO_B O\angle =\pi -\ga -\eta -\frac{\be -\al}{2} =
\pi -\left( \ga +\eta +\frac{\be -\al}{2} \right) .$$
Therefore, by the sine law applied to the triangle $\De _{OO_AO_B }$ we obtain
\begin{equation}\label{e7}
\frac{\sin \left( \ga +\eta +\frac{\be -\al}{2} \right) }{\sin \ga}=
\frac{\sin (\al /2)}{\sin (\be /2)} .
\end{equation}
It is easy to check that $5\cd \sin (1/5)>15/16.$ Since
the function $(\sin x)/x$ is decreasing in $(0,\pi /2)$, it follows that
$(\sin x)/x>15/16$ for every $0<x<1/5$. 
In particular, we have $\sin \be /2 >15\be /32$, and thus 
the right hand side of \eqref{e7} is less than 
$$\frac{\al /2}{15\be /32}=\frac{\al}{\be} \cd \frac{16}{15} <
\frac{3}{4} \cd \frac{16}{15} =\frac{4}{5}.$$
The left hand side of \eqref{e7} equals
$$\cos \left( \eta +\frac{\be -\al}{2} \right) 
+ \cot \ga \cd \sin \left( \eta +\frac{\be -\al}{2} \right) .$$
Since 
$$\cos \left( \eta +\frac{\be -\al}{2} \right) >\cos 2\eta \ge \cos 2/5>9/10,$$
it follows that
$$\cot \ga \cd \sin \left( \eta +\frac{\be -\al}{2} \right) <
 \frac{4}{5} -\frac{9}{10} =-\frac{1}{10}.$$
This implies $\cot \ga <0$, $\pi /2 <\ga <\pi$, and
$$\cot (\pi -\ga ) \cd \sin \left( \eta +\frac{\be -\al}{2} \right) >
 \frac{1}{10}.$$
Since $\sin (\eta +(\be -\al )/2)<\sin 2\eta <2\eta$, we obtain
$\tan (\pi -\ga )<20 \eta.$

The intersection $K_A \cap K_B$ consists of two points, each on
the perpendicular bisector of the segment $O_AO_B $. Let $P$ be the one which
is on the same side of the line $O_AO_B $ as the points $A$ and $B$.
We have
$$O_A O B \angle = O_A O O_B \angle +O_B O B \angle
= \eta +\frac{\be -\al}{2} +\frac{\pi -\be}{2}=
\frac{\pi}{2} +\eta -\frac{\al}{2}.$$
Then, by $\al <\eta$, we obtain that $O_A O B \angle$ is an obtuse angle.
Therefore, we have $|BO_A|>|BO|=1=|BO_B |$.
Since $|OO_A | <|OO_B |$, it follows that 
the perpendicular bisector of the segment $O_AO_B $
intersects the segment $OB$ at a point $Q$. Clearly,
$|O_BQ| < |O_BB| = 1 = |O_BP|$, and thus $P$ and $O_B$ are separated by the line $OB$.
Therefore, $O_A OP \angle >O_A OB\angle  >\pi /2$, and then
$|OP|<|O_A P|=1$. This proves that $P$ is inside $K$.

Since $|AO_A|=|PO_A|=1$, the triangle $\De_{APO_A}$ is isosceles.
If $PO_A A\angle =\de$, then 
$$\ga =OO_A O_B \angle = (PO_A O_B \angle )  +(OO_A A\angle )
-\de < 2\cd (\pi /2)-\de = \pi -\de ,$$
and thus $\de <\pi -\ga$. Therefore, we have
$$|AP|=2\sin (\de /2)<\tan \de<\tan (\pi -\ga)<20 \eta ,$$
which completes the proof of the lemma.
\hfill $\square$

\begin{figure}
\centering
\includegraphics[width=3in]{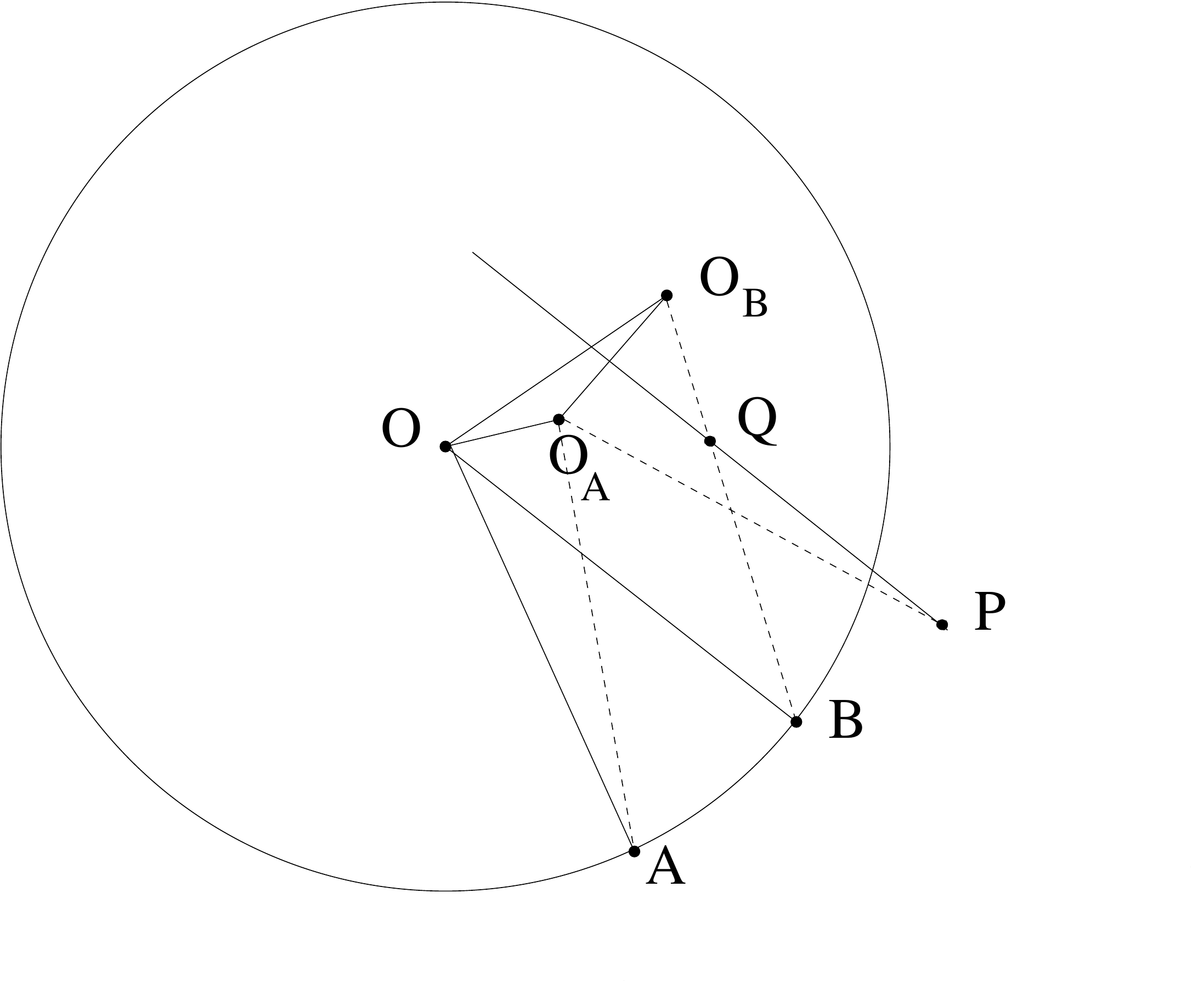}
\caption{Lemma \ref{l9}}
\label{fig8}
\end{figure}

\bsk \noi
{\bf Proof of Lemma \ref{l9}.}
The proof is similar to that of Lemma \ref{l3}
(see Figure \ref{fig8}). 
It is easy to check that $O_AOO_B \angle 
=\eta -(\be -\al )/2$, and that the orientation of
the triangle $OO_A O_B$ is positive. Let $OO_AO_B \angle =\ga$.
Then
$$O_AO_B O\angle =
\pi -\left( \ga +\eta -\frac{\be -\al}{2} \right) ,$$
and from this we obtain $\pi /2 <\ga <\pi$ and
$\tan (\pi -\ga )<20 \eta$ the same way as in the proof of Lemma \ref{l3}.

Since $1=|BO_B |>|BO_A |$, the perpendicular bisector of $\ol{O_A O_B}$
intersects the segment $\ol{O_B B}$ at a point $Q$. Then it follows from
$|O_B Q|<|O_B B|=1$ that $P$ and $O_A$ are separated by the line $\ol{O_B B}$.

We have 
$$OO_A P\angle =2\pi -\left( \ga + (PO_A O_B \angle )  \right) >
2\pi -\left( \pi +\frac{\pi}{2}  \right)= \frac{\pi}{2},$$
and thus $|PO|>|PO_A |=1$. This proves that $P$ is outside $K$.

Since $|O_A P|=|O_B P|=1$, the triangle $\De_{O_A PO_B}$ is isosceles.
Let $\zeta$ denote the angle $O_A PO_B \angle$. 
Since the angle $\ga$ is obtuse, we have 
$|O_A O_B |<|OO_B |$, and thus $\zeta <\be$.

Let $PO_A A\angle =\de$. Then
$$\de =2\pi -\ga  -\frac{\pi -\al}{2} - \frac{\pi - \zeta}{2}
<\pi -\ga +\eta .$$
Hence
\begin{align*}
|AP| & =2\sin (\de /2)<\tan \de<\tan (\pi -\ga +\eta )
=\frac{\tan (\pi -\ga )+\tan \eta}{1-\tan (\pi -\ga )\cd \tan \eta}
< \\
& <\frac{20 \eta +2\eta}{1-40 \eta ^2}<50\eta ,
\end{align*}
which completes the proof of the lemma.
\hfill $\square$

\bsk \noi
{\bf Proof of Lemma \ref{l6}.}
If $\al =2\pi /n$, then the statement follows from the fact that
the disc of radius $|AB|$ can be decomposed into $n$
congruent copies of $H$. Therefore, the statement is also true if $\al$
is a rational multiple of $\pi$. The general case follows by
approximating $\al /\pi$ from below and from above by rational numbers.
\hfill $\square$

\begin{lemma}\label{l14}
Let $a\le b\le c$ be the sides of a triangle, and let $\al$ be its angle
opposite to the side $a$. Then $\al <2a/c.$
\end{lemma}

\proof
Since $\al$ is the smallest angle of the triangle, we have $\al <\pi /2$,
and thus $\al <2\sin \al$. If $\ga$ is the angle opposite to the side $c$,
then we have $\sin \al /\sin \ga =a/c$, $\sin \al \le a/c$ and $\al <2a/c$.
\hfill $\square$

\bsk
In the sequel we shall use the 
following notation. Let $h , \ep$ be positive numbers satisfying
$\ep <10^{-6}$ and $h\le \ep /10^{3}$. 
Let $K_0$ and $K_1$ be the circles of radius $1$ and 
centres $(-\sin (h /2),0)$ and 
$(\sin (h/2),0)$, respectively. The point with coordinates
$(0, \cos (h/2))$ will be denoted by $M$. Thus $M\in K_0 \cap K_1$.

Let $\ol K _i$ denote the disc bounded by $K_i$ $(i=0,1)$.
The closures of the sets 
$\ol K _0 \se \ol K _1$ and $\ol K _1 \se \ol K _0$ are denoted by
$L_0$ and $L_1$. 

The circle of radius $\ep$ and centre $M$ intersects
the lune $L_0$ in the arc $\iv{P_0  P_1}$, where $P_0 \in K_0$ and 
$P_1 \in K_1$.

\bsk \noi
{\bf Proof of Lemma \ref{l4}.}
First we show that 
\begin{equation}\label{e9}
|(P_0 MO\angle )-(\pi /2)|<h+\ep .
\end{equation}
The triangle $\De _{P_0 O_0 M}$ is isosceles, and $|P_0 M|=\ep$.
If $\al _0 = P_0 O_0 M\angle$, then $\al _0 <2\ep$,
by Lemma \ref{l14}.
Since $P_0 MO\angle =((\pi -\al _0 )/2)  +(h/2)$, \eqref{e9} follows.
A similar argument gives
\begin{equation}\label{e10}
|P_1 MO\angle -(\pi /2)|<h+\ep .
\end{equation}
The inequalities \eqref{e9} and \eqref{e10} imply that the angle between
the $x$ axis and any line connecting $M$ to a point of $\iv{P_0  P_1}$ 
is less than $h+\ep$. Therefore, 
the angle between the $y$ axis and any tangent line of $\iv{P_0  P_1}$
is less than $h+\ep$. Now, if a 
line intersects the arc $\iv{P_0  P_1}$ at two points, then it 
is parallel to a tangent line of 
$\iv{P_0  P_1}$, and then the statement of the lemma follows.
\hfill $\square$

\bsk \noi
{\bf Proof of Lemma \ref{l5}.}
First we show that the orientation of the triangle $\De _{Q'' O'C}$ is positive.
Let $D$ denote the circle of centre $O_0$ and radius $1.5\ep$.
Since $|O'O|<\ep$ and $|OO_0 |<h/2$, we have $O' \in D$. Let $\ell$ denote
the line going through $Q''$, tangent to the disc $D$
at a point $E$, and such that 
the orientation of the triangle $\De _{Q'' O_0 E}$ is positive.
Let $M'\in K_0$ be such that the line $\ol{MM'}$
is parallel to $\ell$.
The line $\ol{MM'}$ is obtained from $\ell$ by translating it by a vector
of length $<\ep$. Since the distance between the point $O_0$
and $\ell$ is $1.5\ep$, the distance between the point $O_0$
and $\ol{MM'}$ is at most $2.5\ep$. 

\begin{figure}
\centering
\includegraphics[width=3in]{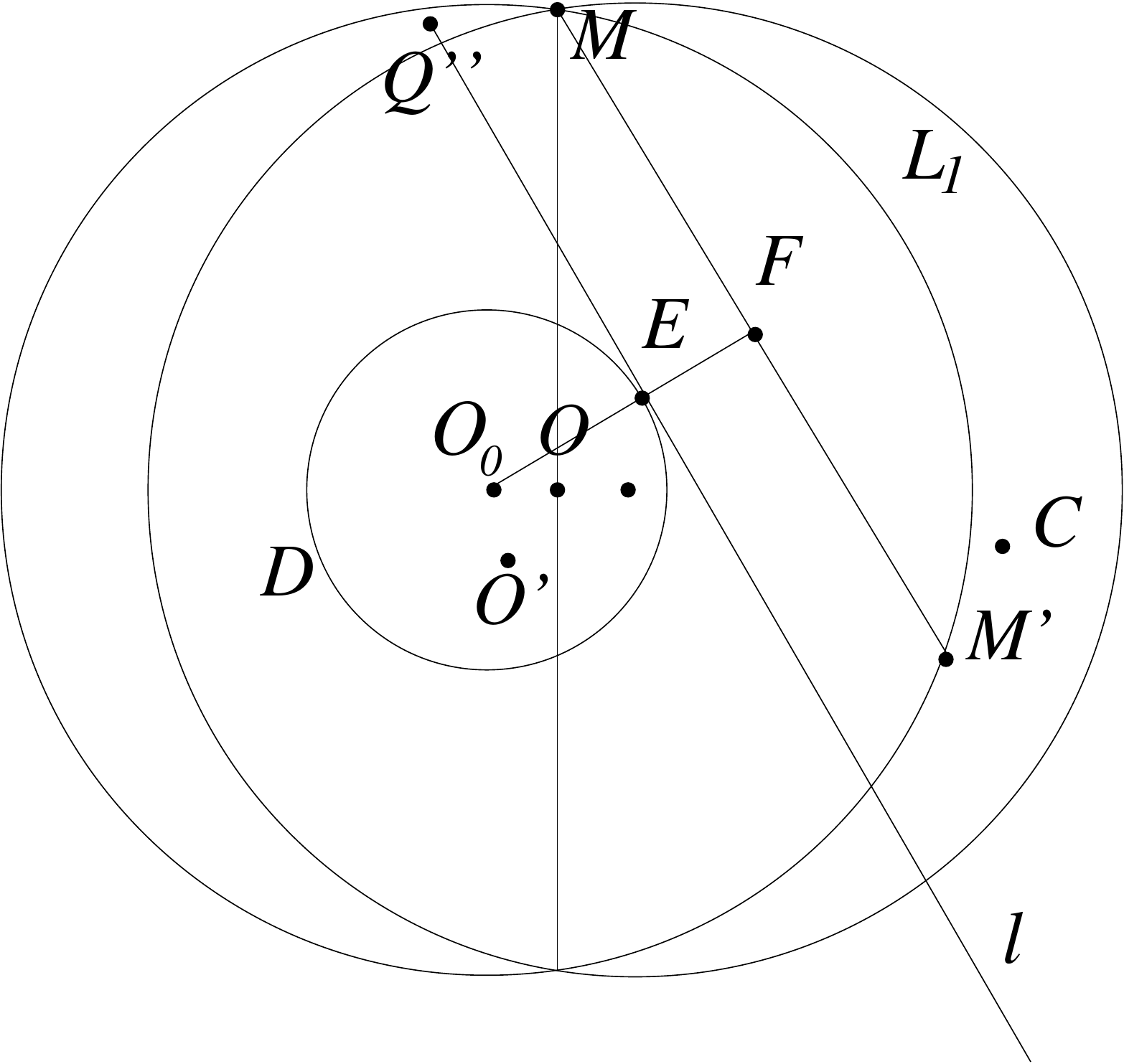}
\caption{Orientation of $\De_{Q'' O'C}$}
\label{fig5}
\end{figure}

Let $F$ denote the middle point of the segment
$MM'$. Then $|O_0 F|\le 2.5\ep$, $|MF|>1-2.5\ep$ and 
$|MM'|=2|MF|>2-5\ep$.
Since $C\in L_1$ and $|MC|\le 2-5\ep$, it follows that
the orientation of the triangle $\De _{MM'C}$ is positive
(see Figure \ref{fig5}). Thus  
the orientation of the triangle $\De _{Q'' O'C}$ is also positive.

The point $Q''$ is inside the circle $K'$, and
the distance $d$ between the point $Q''$ and the circle 
$K'$ is at most the length
of the arc $\iv{Q' Q''}$, which is at most $h\ep$.
Since the distance $d'$ between $C$ and $Q''$
is at least $\ep > 2 \sqrt{h\ep} > 2\sqrt d$,
we may apply
(i) of Lemma \ref{l1}, and find that $\sin \al <2\sqrt{d}< \ep /2$.
Thus $\al <\ep$.

Let the point $\tilde Q$ be selected such that $Q',O',O,\tilde Q$
are the vertices of a parallelogram. Then
$$|\tilde Q M|\le |\tilde Q Q'| +|Q'M| =|O' O| +|Q'M| <2\ep .$$
Thus, by Lemma \ref{l14} applied to the triangle $\De_{\tilde{Q}OM}$ we find that
$\tilde{Q} OM\angle <4\ep$; that is, 
the angle $\ga$ of the line $Q' O'$ and the $y$ axis is less than $4\ep$.

Let $t=|Q'' C|=|\rho (Q'')C|.$ Then $t<|Q'' M|+|MC| <\ep +2-5\ep <2$, and thus
\begin{equation}\label{e31}
|Q''\rho (Q'')|=t\cd 2\sin (\al /2) < 2\cd 2\sin (\al /2) <2\al <2\ep .
\end{equation} 
Therefore, $|Q' \rho (Q'')|\le |Q' Q''|+|Q''\rho (Q'')|<3\ep$, and then
$Q'O'\rho (Q'') \angle <3\ep$.
If the angle between the $x$ axis and the line $\ol{Q' \rho (Q'')}$ is $\zeta$,
then
$$\zeta \le \ga +\left| (O'Q'\rho (Q'') \angle )-(\pi /2)
\right| <4\ep + \tfrac{1}{2} \cd Q'O'\rho (Q'') \angle  <6\ep ,$$
which proves (iii).\hfill $\square$

\bsk \noi
\begin{lemma} \label{l7}
There exists an $n_0$ depending on $h$
but not on $n$ such that
$|XY| < 6h$ for every $n\ge n_0$, $0\le i< n$, $X\in \iv{A_i^0 A_i^1}$
and $Y\in \iv{A_{i+1}^0 A_{i+1}^1}$.
\end{lemma}

\proof
Since the function $\arccos x$ is uniformly continuous on the interval
$\nl$, there is a $u>0$ such that $|\arccos x -\arccos y|<h$
whenever $x,y\in \nl$ and $|x-y|\le u$. We show that $n_0 =2/u$
satisfies the requirement. Suppose $n\ge n_0$.

For every $X\in \iv{A_i^0 A_i^1}$ the distance $|XA_i^0 |$ is less than
the length of the arc $\iv{A_i^0 A_i^1}$ which is less than $2h$.
Similarly, we have $|YA_{i+1}^0 |<2h$ for every $Y\in \iv{A_{i+1}^0 A_{i+1}^1}$.

Let $\al _i , \al _{i+1},\de _i$ denote the angles $O_0 MA_{i}^0 \angle$,
$O_0 MA_{i+1}^0 \angle$, $A_{i+1}^0 MA_i^0 \angle$, respectively. Then
\begin{equation}\label{e12}
\de _i  =\al _i -\al _{i+1}=\arccos (r_i /2) -\arccos (r_{i+1} /2)<h,
\end{equation} 
as $|r_{i+1} -r_i | =R/n<2/n<u$. Since $A_i^0 O_0 A_{i+1}^0 \angle =2\de _i$,
\eqref{e12} implies
$|A_i^0 A_{i+1}^0| =2\sin \de _i  <2\de _i <2h$, and hence,
$$|XY|\le |XA_i^0 | +|A_i^0 A_{i+1}^0| +|YA_{i+1}^0 |<6h.$$
\hfill $\square$

\bsk \noi
{\bf Proof of Lemma \ref{l8}.}
Let $n_0$ be as in Lemma \ref{l7}. We show that 
$c=\max(2n_0 ,8/\ep )$ satisfies the requirement.
If $n\le n_0$ then $|C_i^x C_{i+1}^x |\le 2\le c/n$
and $|C_i^x C_{i+1}^{x+2^{-i-1}} |\le 2\le c/n$ 
for every $i$ and $x$. Therefore, we may assume that $n>n_0$.

We put $|O_0^x C_i^x |=a$, $|C_i^x M |=b$
and  $|MO_0^x |=g$. 
Note that $a=1$ by $C_i^x \in K_0^x$, and $b=r_i \le 2-5\ep$. 
Let
$C_{i}^x M O_0^x \angle =\al$, $MO_0^x C_{i}^x  \angle =\be$, 
$O_0^x C_{i}^x M \angle =\ga$. 
Our first aim is to show $\ga >\ep$.

Since
$|O_0^x O|<\ep$ and $|OO_1 |<h/2$, we have $|O_0^x O_1 |<2\ep$
and $g=|O_0^x M |>1-2\ep$. Since $|O_1 M|=1$, we find
$O_0^x MO_1 \angle <4\ep$ by Lemma \ref{l14}.
We also have $O_1 MC_i^x \angle < \pi /2$, since $\ol{O_1 M}$ is a
ray of $K_1$ and $C_i^x$ is either inside or on the circle $K_1$. Therefore, 
we obtain $\al <(\pi /2)+4\ep$.
If $b<1-2\ep$, then $b$ is the shortest side of the triangle
$\De _{O_0^x C_{i}^x M}$. Then, in this case, $\be <\pi /3$ and
$\ga =\pi -\al -\be > (\pi /6)-4\ep >\ep$.

Now suppose $b\ge 1-2\ep$.
Then we have, by $a=1$ and $b\le 2-5\ep$,
\begin{equation*}
\begin{split}
1-\cos \ga &=1-\frac{a^2 +b^2 -g^2}{2ab}=\frac{g^2 -(a-b)^2}{2ab}
\ge \frac{(1-2\ep )^2 -(1-5\ep )^2}{2ab}=\\
&\ge  \frac{(3\ep)(2-7\ep )}{4}>\ep .
\end{split}
\end{equation*} 
Therefore, we obtain $\ga >1-\cos \ga >\ep$.
\begin{figure}
\centering
\includegraphics[width=3.5in]{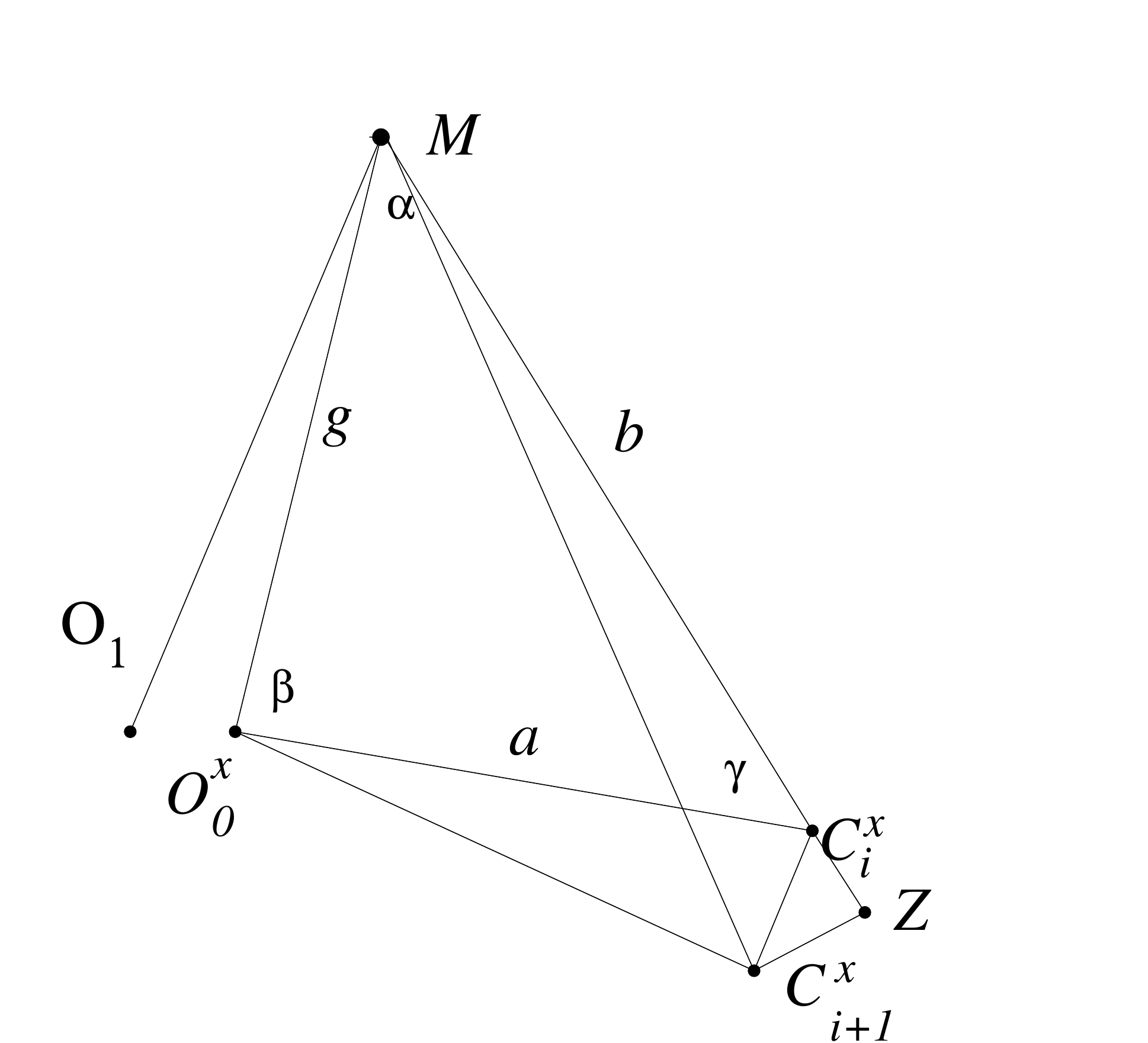}
\caption{Lemma \ref{l8}}
\label{fig6}
\end{figure}

Let $Z\in \iv{A_{i+1}^0 A_{i+1}^1}$ be such that $M$, $C_i^x$
and $Z$ are collinear. Note that $|C_{i+1}^xC_{i}^x | <6h$
by Lemma \ref{l7} and $n>n_0$.
Since the triangle $\De _{O_0^x C_{i+1}^x C_{i}^x}$ 
is isosceles and $|O_0^x C_{i+1}^x |=|O_0^x C_{i}^x |=1$,
we have $C_{i+1}^x O_0^x C_{i}^x \angle <12h$ by Lemma \ref{l14}.
Then $O_0^x C_{i}^x C_{i+1}^x  \angle >(\pi /2)-6h$, and
$$C_{i+1}^x C_{i}^x Z  \angle <\pi - ((\pi /2)-6h)-\ep <(\pi  /2)- (\ep /2).$$
Since $C_{i+1}^x Z C_{i}^x \angle <\pi /2$, we obtain
$C_{i}^x  C_{i+1}^x Z\angle >\ep /2$. Therefore, 
\begin{align*}
\frac{2}{n} >\frac{R}{n} &=| C_{i}^x Z|\ge |C_{i+1}^x C_i^x |\cd 
\sin (C_{i}^x  C_{i+1}^x Z\angle )>|C_{i+1}^x C_i^x |\cd \sin (\ep /2)  >\\
&> |C_{i+1}^x C_i^x |\cd \frac{\ep}{4} .
\end{align*}
This proves $|C_{i+1}^x C_i^x |<8/(\ep n) \le c/n$.

The inequality
$|C_i^x C_{i+1}^{x+2^{-i-1}} |\le c/n$ can be proved similarly;
we have to replace the point $O_0^x$ by $O_1^{x+2^{-i}}$,
and use the fact that $K_1^{x+2^{-i}}$ contains the points 
$C_i^x$ and $C_{i+1}^{x+2^{-i-1}}$.
\hfill $\square$

\begin{figure}
\centering
\includegraphics[width=4in]{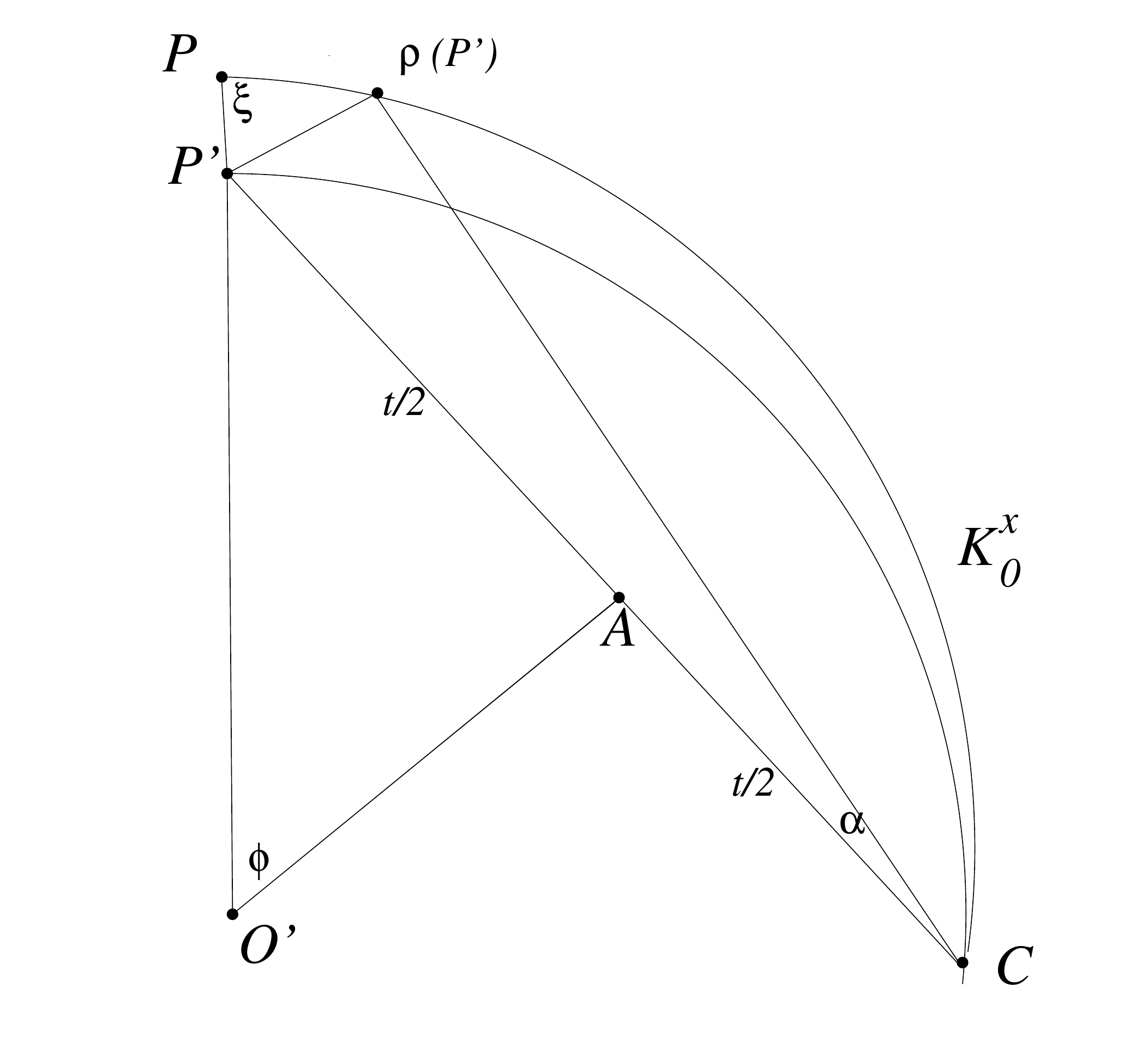}
\caption{Lemma \ref{l10}}
\label{fig11}
\end{figure}

\bsk \noi
{\bf Proof of Lemma \ref{l10}.}
We shall use the abbreviations $P=P_x$, $P'=P_{x+2^{-i}}$,
$O'=O_1^{x+2^{-i}}$, $C=C_i^x$ and $\al =\al _i^x$ (see Figure \ref{fig11}). 
Also, we denote by $\rho$ the rotation about the point $C$ mapping 
$K_1^{x+2^{-i}}$ onto $K_0^x$.
Since $P'$ and $C$ lie on the circle $K_1^{x+2^{-i}}$ of centre $O'$,
we have $|O'P'|=|O'C|=1$, and thus the triangle $\De _{P'O'C}$
is isosceles. Let $A$ denote the middle point of the segment
$\ol{P'C}$. Then triangle $\De _{P'O'A}$ has a right angle at the vertex
$A$.

We have $|P'M|=\ep$ and $|O'O|<\ep$ by 
\eqref{e14}.
Let the point $Q$ be selected such that $O',O,Q,P'$
are the vertices of a parallelogram. Then
$$|Q M|\le | Q P'| +|P'M| =|O' O| +|P'M| <2\ep .$$
Thus, by Lemma \ref{l14} applied to the triangle $\De_{MOQ}$ we find that
the angle between the $y$ axis and the segment $\ol{O'P'}$
is less than $4\ep$. 

We denote the angle
$P'O'A \angle$ by $\phi$. Clearly, $\phi =\arcsin (t/2).$

By Lemma \ref{l4}, 
the angle between the $y$ axis and the segment $\ol{PP'}$
is less than $2\ep$. Therefore, if we add the angles at the vertex $P'$
then we obtain
\begin{equation}\label{e28}
\left| \left( \frac{\pi}{2} -\phi \right) + (\rho (P') P'C \angle )+
(PP'\rho (P') \angle ) -\pi \right| <6\ep .
\end{equation} 
We have $|P'C|=|\rho (P')C|=t$, as $\rho$ is a rotation about the point
$C$. Thus the triangle $\De _{P'C\rho (P')}$ is isosceles, and
$\rho (P') P'C \angle =(\pi -\al )/2$. If we denote the angle
$PP'\rho (P') \angle$ by $\psi$, then \eqref{e28} gives
$|\psi -\phi -(\al /2)|<6\ep$. Now we have $\al <\ep$ by 
(i) of Lemma \ref{l5}, and thus $|\psi -\arcsin (t/2)|<7\ep$.
Therefore, we have $|\cos (\arcsin (t/2))-\cos \psi | <7\ep$; that is, 
\begin{equation}\label{e29}
|\cos \psi - \sqrt{1-(t/2)^2} | <7\ep .
\end{equation}
Since $t\le R + \ep < 2-4\ep$, we have 
$$\sqrt{1-(t/2)^2} >\sqrt{1-(1-2\ep )^2}  =\sqrt{4\ep -4\ep ^2}  >
\sqrt{\ep}  >10^3 \ep$$
by $\ep <10^{-6}$. Then \eqref{e29} gives $\cos \psi  >900\ep$.
We also have
\begin{equation}\label{e37}
0.99  <1-\frac{7}{10^3}< \frac{\cos \psi}{\sqrt{1-(t/2)^2}}<
1+\frac{7}{10^3}<1.01.
\end{equation} 
The angle between the $x$ axis and segment $\ol{P\rho (P')}$
is less than $6\ep$ by (ii) of Lemma \ref{l5}. If $\xi$ denotes
the angle $P'P\rho (P')\angle$, then this implies $|\xi -(\pi /2)|<8\ep$
and $|\cot \xi |<10\ep$ by $\ep <10^{-6}$.

Now we apply the sine law for the triangle $\De _{PP'\rho (P')}$.
Since $|P'\rho (P')|=2t\sin (\al /2)$, we obtain
\begin{equation}\label{e30}
\frac{|PP'|}{2t\sin (\al /2)} =\frac{\sin (\pi-\xi -\psi)}{\sin \xi}=
\frac{\sin (\psi+\xi)}{\sin \xi}=\cos \psi +\cot \xi \cd \sin \psi .
\end{equation} 
By $\cos \psi >900\ep$ and $|\cot \xi |<10 \ep$ we find that
$$0.95 \cos \psi <\cos \psi +\cot \xi \cd \sin \psi <1.02 \cos \psi .$$
The length of the arc $\iv{PP'}$ equals $h\ep \cd 2^{-i}$, and thus
$|PP'|=2\ep \sin (h\cd 2^{-i-1})$. This gives
$$0.95 h\ep \cd 2^{-i} <|PP'| <h\ep \cd 2^{-i}.$$
Therefore, by \eqref{e37} and \eqref{e30} we obtain
\begin{align*}
\al & >2\sin (\al /2) =\frac{|PP'|}{t(\cos \psi +\cot \xi \cd \sin \psi )}
\ge \frac{0.95 \cd h\ep \cd 2^{-i}}{1.02 \cd  t\cd \cos \psi}
>0.93 \frac{h\ep \cd 2^{-i}}{t\cd \cos \psi} >\\
&> 0.9 \frac{h\ep \cd 2^{-i}}{t\cd \sqrt{1-(t/2)^2}} 
\end{align*}
and 
\begin{align*}
\al &< 1.01 \cd 2\sin (\al /2) =\frac{1.01 \cd  |PP'|}{t(\cos \psi +
\cot \xi \cd \sin \psi )} < 
\frac{1.01}{0.95}  \cd \frac{h\ep \cd 2^{-i}}{t\cd \cos \psi}<\\
&<  1.07 \cd \frac{h\ep \cd 2^{-i}}{t\cd \cos \psi}< 
 1.1 \cd\frac{h\ep \cd 2^{-i}}{t\cd \sqrt{1-(t/2)^2}} ,
\end{align*}
which proves \eqref{e22}.

\bsk \noi
{\bf Proof of Lemma \ref{l+}.}
We shall prove that the statement holds with $q=1-\ep^4 $. 

For all $x \in D_n, x >0$, let $C_n^xMC_n^{x-2^{-n}} \angle=(\beta_n^{x})' $. 
Then we have 
$$\sum_{x \in D_n, \ x>0} (\beta_n^{x} )'=A_n^0 M A_n^1 \angle=h.$$ 
Therefore, it is enough to prove that there is a positive constant $ q < 1$ 
depending only on $\ep$ such that for all $x \in D_n, 0 < x <1$, 
$\beta_n^{x} \leq q \cd (\beta_n^{x} )'$. 

 Fix $x \in D_x, 0 < x < 1$. We shall use the abbreviations 
$K=K_0^x, L=K^x_1, P=P^x, A=C_n^x, B=C_n^{x-2^{-n}}, \beta=\beta_n^x$, and 
$\beta'=( \beta_n^{x} )'$. Then we have that $K$ contains $P$ and $A$, 
$L$ contains $P$ and $B$. Let $\rho$ denote the rotation about the point $P$ 
by angle $\beta$ in the positive direction. Then $\rho(L)=K$. 
Put $B'=\rho(B) \in K$, $b=|BB'|$, $s=|BP|=|B'P|$, then $BPB' \angle=\be$ and 
$b=2s \cd \sin (\beta /2)$. 

It is easy to check that $\beta /2  < \ep $, so we have
$$b=2s\sin \left(\frac{\beta}{2} \right) > 
2s \cd  \frac{\beta}{2} \cdot \frac{\sin \ep}{\ep} > 
s \cdot \beta \cdot \left(1-\frac{\ep^2 }{6} \right)$$
and 
$$\beta < \frac{b}{s} \cdot \frac{1}{1- (\ep^2 /6)} < 
\frac{b}{s} \cdot \left( 1+ \ep^2 \right).$$ Thus we have 

\begin{equation}\label{ll}
\frac{\be}{\be'}< \frac{b}{\be'} \cdot \frac{1}{s} \left( 1+ \ep^2 \right).
\end{equation}

We shall estimate $b/\be'$ from above and  $s$ from below. 
First we estimate $s$.

Put $\delta=BMP \angle = NMP \angle + BMN \angle$. 
If the $y$ coordinate of $P$ is greater than or equal to the $y$ 
coordinate of $M$ (which is $\cos(h/2)$), then $NMP \angle >\pi/2$. 
Otherwise the tangent line at $P$ of the circle with centre $M$ 
and radius $\ep$ intersects the 
$y$ axis at a point $Q$ with $y$ coordinate less than $\cos(h/2)$. Then 
$$NMP \angle = QMP \angle =\pi - \pi/2 - PQM \angle  
\geq \pi/2 - h - \ep$$ by lemma \ref{l4}. 

The point $B$ is on the arc $\iv{A_n^0 A_n^1}$ lying in the lune
$L_1$. This implies $BMO_0 \angle \ge A_n^0 MO_0 \angle$. 
Since the triangle $\De_{O_0 A_n^0 M}$ is 
isosceles with $|O_0 M|=|O_0 A_n^0|=1$ and $|MA_n^0|=R$, we have 
$$BMO_0 \angle \ge A_n^0 M O_0 \angle = \frac{\pi}{2}-\arcsin (R/2).$$
Thus 
$$BMN \angle =BMO_0 \angle -NMO_0 \angle \ge 
\frac{\pi}{2}-\arcsin (R/2) -(h/2).$$
Therefore, we have the estimate
$$ \de \geq \pi - \arcsin  (R/2) - (3h/2)-\ep 
\geq \pi - \arcsin (R/2) - 2 \ep$$
by $\ep \geq 2h$. Then
\begin{align*}
\cos \de & \leq \cos \left(\pi - \arcsin (R/2) - 2 \ep \right)=\\
&= -\sqrt{1-(R^2 /4)} \cdot \cos(2\ep) + (R/2)\cd  \sin(2\ep) \le\\
&\le  -\sqrt{1-(R^2 /4)} (1-2\ep^2)+R \ep \leq \\
&\le -\sqrt{1-(R^2 /4)} + 4 \ep
\end{align*}
by $1 \leq R <2$.
Now we apply the cosine law for $\De_{PBM}$:
$$s^2=R^2+\ep^2-2R \ep \cos \de \geq R^2 + \ep^2 + 
2R \ep \sqrt{1-(R^2 /4)} - 8R \ep^2. $$
Using $1 \leq R < 2$ again, we have 
$$s^2 \geq  R^2 + R \ep \sqrt{4-R^2} -15 \ep^2 \geq 
\left( R + (\ep /2) \cd  \sqrt{4-R^2} \right)^2 -16 \ep^2.$$
Thus
$$s- \left(R + (\ep /2) \cd \sqrt{4-R^2} \right) \geq \frac{-16 \ep^2}
{s+R + (\ep /2)\cd \sqrt{4-R^2}}  \geq \frac{-16 \ep^2}{2}=-8\ep^2,$$
where we estimated the denominator by using $s \geq R \geq 1$. We get
\begin{equation}\label{eqs}
s \geq R + (\ep /2) \cd \sqrt{4-R^2} -8\ep^2.
\end{equation}
Now we estimate $b/\be'$ from above. 
Put $a=|AB|, d=|AB'|, t=|PA|$. 
Applying the cosine low for $\De_{PAM}$, we have
$$t^2= R^2+\ep^2-2 R \ep \cos(\de + \be').$$
Thus
\begin{equation}\label{eqt-s}
t-s=\frac{t^2-s^2}{t+s}=\frac{2 R \ep (\cos(\de)-\cos(\de+\be'))}{t+s} 
\leq \frac{2R\ep \be'}{2R}=\ep \be' ,
\end{equation}
where we used $t \geq s \geq R$ and 
the Lipschitz continuity of the function $\cos x$.

Now we estimate $d$. 
The points $P, B',$ and $A$ lie on the circle $K$ of radius $1$, and 
$|PA|=t, |PB'|=s$, thus we get
$$d < \iv{AB'} = \iv{PA}-\iv{PB'}=
2 \arcsin \left(\frac{t}{2} \right) - 2 \arcsin \left(\frac{s}{2} \right).$$ 
Applying the mean value theorem for the function $\arcsin(x)$ and 
using \eqref{eqt-s}, we get 

\begin{equation}\label{eqd}
d < 2 \arcsin \left(\frac{t}{2} \right) - 2 \arcsin \left(\frac{s}{2} \right) \leq 
2 \left(\frac{t}{2}-\frac{s}{2} \right) \cdot \frac{1}{\sqrt{1-\frac{t^2}{4}}} 
\le \frac{2\ep \be'}{\sqrt{4-t^2}}.
\end{equation}

Put
$\phi=BAB' \angle = BAP \angle + PAB' \angle.$
Since $PAB'\angle $ is smaller than the angle $\psi$ 
of the secant $PA$ and the tangent line 
of the circle $K$ at $A$, and $\psi= \arcsin (t/2)$,
we have that $PAB' \angle \leq \arcsin (t/2).$ 
On the other hand, one can easily check that $BAP \angle < BAM \angle$, 
and since $\De_{BAM}$ is isosceles, $BAM \angle < \pi/2$. 
Thus $\phi  \leq (\pi/2) + \arcsin (t/2)$ 
and  
$$\cos\phi \geq \cos  \left( (\pi/2) + \arcsin (t/2) \right)=-t/2.$$
Now we shall apply the cosine law for $\De_{ABB'}$. 
$$b^2=a^2+d^2-2ad\cos(\phi) \leq a^2+d^2+adt.$$ 
Since $BA$ is a chord of the circle with centre $M$ and 
radius $R$ having central angle $\be'$, 
we have that $a=|BA|< \iv{BA}=R \beta'$. 
Using this and \eqref{eqd}, we get the estimates
$$b^2 < (R \beta')^2+ \frac{4(\ep \be')^2}{4-t^2}+
\frac{2Rt\ep (\be')^2}{\sqrt{4-t^2}}$$
and
$$\frac{b^2}{\be'^2} < R^2 + 
\frac{4\ep^2 }{4-t^2}+\frac{2Rt\ep }{\sqrt{4-t^2}}=
\left(R+ \frac{t\ep}{\sqrt{4-t^2}}\right)^2 + \ep^2.$$
Now we are ready to estimate $b/\be'$.
$$\frac{b}{\be'}-\left(R+\frac{t\ep}{\sqrt{4-t^2}}\right)=
\frac{  (b/\be ')^2 - \left(R+\frac{t\ep}{\sqrt{4-t^2}}\right)^2}
{ (b/\be ') +R+\frac{t\ep}{\sqrt{4-t^2}}}<\frac{\ep^2}{1}=\ep^2,$$
where we used that $R \geq 1$. Thus
\begin{equation}\label{eqb/be'}
\frac{b}{\be'}<R+\frac{t\ep}{\sqrt{4-t^2}} +\ep^2 .
\end{equation}
Now we estimate the quotient $\frac{\beta}{\beta'}$ 
using \eqref{ll}, \eqref{eqs} and \eqref{eqb/be'}. 
\begin{align*}
\frac{\be}{\be'}& < \frac{b}{\be'} \frac{1}{s} \left ( 1+ \ep^2 \right)<
\frac{R+\frac{t\ep}{\sqrt{4-t^2}} +\ep^2 }
{R + (\ep /2) \cd \sqrt{4-R^2} \ep - 8 \ep^2 }
\left ( 1+ \ep^2  \right)= \\
&=\left(1- \frac{ (\ep /2) \cd \sqrt{4-R^2} -\frac{t\ep}{\sqrt{4-t^2}} -9 \ep^2}
{R + (\ep /2) \cd \sqrt{4-R^2} -8 \ep^2} \right) \left ( 1+ \ep^2  \right).
\end{align*}

If we show that 
$$\frac{ (\ep /2) \cd \sqrt{4-R^2} -\frac{t\ep}{\sqrt{4-t^2}} -9 \ep^2}
{R + (\ep /2) \cd \sqrt{4-R^2}-8 \ep^2} \geq \ep^2,$$ then we are done since then 
$$\frac{\beta}{\beta'} < (1-\ep^2) ( 1+ \ep^2) = 
1-\ep^4=q < 1.$$
Since $R + (\ep /2) \cd \sqrt{4-R^2} -8 \ep^2 \leq 2+ 1=3$, 
it is enough to show that 
$$ (\ep /2) \cd \sqrt{4-R^2} -\frac{t\ep}{\sqrt{4-t^2}} -9 \ep^2 \geq 3 \ep^2.$$
In turn, by $\sqrt{4-R^2} \geq \sqrt{4-t^2}$, it is enough to show that 
$$(\ep /2) \cd \sqrt{4-t^2} -\frac{t\ep}{\sqrt{4-t^2}} -9 \ep^2 \geq 3 \ep^2, $$ 
which is equivalent to
\begin{equation}\label{e100}
\frac{4-t^2-2t}{\sqrt{4-t^2}} \geq 24 \ep.
\end{equation}
Since $R= 1.227$, we have $t < R + \ep < 1.227 + 10^{-6}< 1.228$.
Then $t<1.23<\sqrt 5 -1$, and thus $4-t^2-2t>0$.
Therefore, \eqref{e100} is equivalent to 
$$(4-t^2-2t)^2 \geq (4-t^2) (24 \ep)^2.$$
In order to check this inequality, it is enough to show that 
\begin{equation}\label{eqend}
t^4+4t^3-4t^2-16t+15.999 \geq 0,
\end{equation}  
considering that that $\ep < 10^{-6}$. 
One can check that if $0 < t \leq 1.228$, then \eqref{eqend} is indeed 
satisfied, and we are done. 
\hfill $\square$

\vfill \eject

\bsk 
\noi
{\small Department of Analysis, Institute of Mathematics,
E\"otv\"os Lor\'and University

Budapest, P\'azm\'any P\'eter s\'et\'any 1/C, 1117 Hungary

{\rm e-mail:} herakornelia@gmail.com (Korn\'elia H\'era)

laczk@cs.elte.hu (Mikl\'os Laczkovich)}

\end{document}